\documentclass{gtart}

\def\ifplaintex{\expandafter\ifx\csname documentclass\endcsname\relax}

\def\gtp{{\mathsurround=0pt\it $\cal G\mskip-2mu$eometry \&\ 
$\cal T\!\!$opology $\cal P\!$ublications}}  

\def\recd{{\small Received:\qua\receiveddate\ifx\reviseddate\relax
\else\qquad Revised:\qua\reviseddate\fi\par}} 

\def\lognumber#1{\def\thelognumber{#1}}
\def\volumenumber#1{\def\thevolumenumber{#1}}
\def\volumeyear#1{\def\thevolumeyear{#1}}
\def\papernumber#1{\def\thepapernumber{#1}}
\def\pagenumbers#1#2{\def\startpage{#1}\def\finishpage{#2}}
\def\published#1{\def\publishdate{#1}}

\def\received#1{\def\receiveddate{#1}}
\def\revised#1{\def\reviseddate{#1}}
\def\accepted#1{\def\accepteddate{#1}}
\def\asciititle#1{\def\theasciititle{#1}}

\long\def\asciiabstract#1{\long\def\theasciiabstract{#1}}

\let\\\par\let\thelognumber\relax\let\thevolumenumber\relax
\let\thepapernumber\relax\let\thevolumeyear\relax\let\startpage\relax
\let\finishpage\relax\let\publishdate\relax\let\receiveddate\relax
\let\reviseddate\relax\let\accepteddate\relax\let\theasciititle\relax
\let\theasciiauthors\relax
\let\theasciiabstract\relax

\let\theasciiemail\relax

\ifplaintex
\font\logobig=cmssbx10 scaled 3836
\font\logomed=cmssbx10 scaled 2557
\else
\font\logobig=cmssbx10 scaled 4200
\font\logomed=cmssbx10 scaled 2800
\fi

\long\def\makeagttitle{   
\count0=\startpage
\agt\hfill      
\hbox to 45truept{\vbox to 0pt{\vglue -13truept{\logomed A\kern -.37em{\logobig 
T}\kern -.38em G}\vss}\hss}
\break
{\small Volume \thevolumenumber\ (\thevolumeyear)
\startpage--\finishpage\nl
Published: \publishdate}

\vglue .25truein

{\parskip=0pt\leftskip 0pt plus
1fil\def\\{\par\smallskip}{\Large\bf\thetitle}\par\medskip} \vglue
0.05truein

{\parskip=0pt\leftskip 0pt plus 1fil\def\\{\par}{\sc\theauthors}
\par\medskip}%
 
\vglue 0.03truein 

{\small\leftskip 25truept\rightskip 25truept{\bf Abstract}\stdspace\theabstract

{\bf AMS Classification}\stdspace\theprimaryclass
\ifx\thesecondaryclass\relax\else; \thesecondaryclass\fi\par
{\bf Keywords}\stdspace \thekeywords\par}\vglue 7truept

}   

\ifplaintex
\font\phead=cmsl9 scaled 950
\font\pnum=cmbx10 scaled 913
\font\pfoot=cmsl9 scaled 950
\headline{\vbox to 0pt{\vskip -4.5mm\line{\small\phead\ifnum
\count0=\startpage ISSN 1472-2739 (on-line) 1472-2747 (printed)
\hfill {\pnum\folio}\else\ifodd\count0\def\\{ }%
\ifx\theshorttitle\relax\thetitle\else\theshorttitle\fi\hfill{\pnum\folio}
\else\def\\{ and }{\pnum\folio}\hfill\ifx\theshortauthors\relax\theauthors
\else\theshortauthors\fi\fi\fi}\vss}}
\footline{\vbox to 0pt{\vglue 0mm\line{\small\pfoot\ifnum\count0=\startpage
\copyright\ \gtp\hfill\else
\agt, Volume \thevolumenumber\ (\thevolumeyear)\hfill\fi}\vss}}
\else
\font=cmsl9 scaled 1050
\font\lnum=cmbx10 
\font=cmsl9 scaled 1050
\makeatletter
\def\@oddhead{{\small\lhead\ifnum\count0=\startpage ISSN 1472-2739 
(on-line) 1472-2747 (printed)\hfill {\lnum\number\count0}\else\ifodd\count0
\def\\{ }\ifx\theshorttitle\relax \thetitle \else\theshorttitle\fi\hfill
{\lnum\number\count0}\else\def\\{ and }{\lnum\number\count0}
\hfill\ifx\theshortauthors\relax 
\theauthors\else\theshortauthors\fi\fi\fi}}\def\@evenhead{\@oddhead}
\def\@oddfoot{\small\lfoot\ifnum\count0=\startpage\copyright\ \gtp\hfill\else
\agt, Volume \thevolumenumber\ (\thevolumeyear)\hfill\fi}
\def\@evenfoot{\@oddfoot}
\makeatother
\fi
\let\maketitlepage\makeagttitle
\let\makeshorttitle\maketitlepage
\let\maketitle\maketitlepage

\newwrite\gtoutfile
\long\gdef\makeheadfile{  
{\def\\{, }\def\s{ }
\immediate\openout\gtoutfile head.xxx
\immediate\write\gtoutfile{To: math@arxiv.org}
\immediate\write\gtoutfile{Subject: put OR rep NNNNN:ppppp}
\immediate\write\gtoutfile{--text follows this line--}
\immediate\write\gtoutfile{Proxy-for: \ifx\theasciiauthors\relax
\theauthors\else\theasciiauthors\fi\s<\ifx\theasciiemail\relax\theemail\else\theasciiemail\fi>}
\immediate\write\gtoutfile{\noexpand\\}
\immediate\write\gtoutfile{Authors: \ifx\theasciiauthors\relax
\theauthors\else\theasciiauthors\fi}
{\def\\{ }\immediate\write\gtoutfile{Title: \ifx\theasciititle\relax
\thetitle\else\theasciititle\fi}}
\immediate\write\gtoutfile{Subj-class: GT or SG, GR etc}
\immediate\write\gtoutfile{MSC-class: \theprimaryclass\ifx\thesecondaryclass\relax\else, \thesecondaryclass\fi}
\immediate\write\gtoutfile{Journal-ref: Algebr. Geom. Topol. \thevolumenumber\s
(\thevolumeyear) \startpage-\finishpage}
\immediate\write\gtoutfile{Comments: Published by Algebraic and
Geometric Topology at}
\immediate\write\gtoutfile{\s\s\s  http://www.maths.warwick.ac.uk/agt/AGTVol\thevolumenumber/agt-\thevolumenumber-\thepapernumber.abs.html}
\immediate\write\gtoutfile{\noexpand\\}
\immediate\write\gtoutfile{}
\ifx\theasciiabstract\relax
\immediate\write\gtoutfile{\theabstract}\else
\immediate\write\gtoutfile{\theasciiabstract}\fi
\immediate\write\gtoutfile{}
\immediate\write\gtoutfile{\noexpand\\}
\immediate\write\gtoutfile{}
\immediate\closeout\gtoutfile}}  

\def\maketitlepage{\makeagttitle\makeheadfile}
\let\makeshorttitle\maketitlepage
\let\maketitle\maketitlepage

\lognumber{31}
\volumenumber{3}
\volumeyear{2003}
\papernumber{31}
\published{2 October 2003}
\pagenumbers{921}{968}
\received{26 February 2003}
\revised{17 July and 2 September 2003}
\accepted{5 September 2003}

\usepackage{amsmath,amssymb}

\usepackage{graphicx}

\theoremstyle{plain}

\newtheorem{cor}{Corollary}[subsection]
\newtheorem{lem}[cor]{Lemma}

\newtheorem{prop}[cor]{Proposition}
\theoremstyle{definition}
\newtheorem{defi}[cor]{Definition}
\newtheorem{rem}[cor]{Remark}

\theoremstyle{plain}
\newtheorem{thm}{Theorem}
\theoremstyle{definition}
\newtheorem{example}[cor]{Example}

\newcommand{\sign}{\operatorname{sign}}

\newcommand{\lat}{\operatorname{lat}}

\newcommand{\Z}{\mathbb{Z}}

\newcommand{\s}{\mathcal{S}}

\def\theenumi{\roman{enumi}}

\begin{document}
\title{Algebraic linking numbers of knots in 3--manifolds}
\asciititle{Algebraic linking numbers of knots in 3-manifolds}
\author{Rob Schneiderman}
\address{Courant Institute of Mathematical Sciences\\New York 
University\\251 Mercer Street\\New York NY 10012-1185, USA}
\email{schneiderman@courant.nyu.edu}

\begin{abstract}
Relative self-linking and linking ``numbers'' for pairs of oriented
knots and 2--component links in oriented 3--manifolds are defined in
terms of intersection invariants of immersed surfaces in
4--manifolds. The resulting concordance invariants generalize the
usual homological notion of linking by taking into account the
fundamental group of the ambient manifold and often map onto
infinitely generated groups. The knot invariants generalize the type 1
invariants of Kirk and Livingston and when taken with respect to
certain preferred knots, called {\em spherical knots}, relative
self-linking numbers are characterized geometrically as the complete
obstruction to the existence of a singular concordance which has all
singularities paired by Whitney disks. This geometric equivalence
relation, called {\em $W\!$--equivalence}, is also related to finite
type 1--equivalence (in the sense of Habiro and Goussarov) via the
work of Conant and Teichner and represents a ``first order''
improvement to an arbitrary singular concordance. For null-homotopic
knots, a slightly weaker equivalence relation is shown to admit a
group structure.
\end{abstract}

\asciiabstract{Relative self-linking and linking `numbers' for pairs
of oriented knots and 2-component links in oriented 3-manifolds are
defined in terms of intersection invariants of immersed surfaces in
4-manifolds. The resulting concordance invariants generalize the usual
homological notion of linking by taking into account the fundamental
group of the ambient manifold and often map onto infinitely generated
groups. The knot invariants generalize the type 1 invariants of Kirk
and Livingston and when taken with respect to certain preferred knots,
called spherical knots, relative self-linking numbers are
characterized geometrically as the complete obstruction to the
existence of a singular concordance which has all singularities paired
by Whitney disks. This geometric equivalence relation, called
W-equivalence, is also related to finite type 1-equivalence (in the
sense of Habiro and Goussarov) via the work of Conant and Teichner and
represents a `first order' improvement to an arbitrary singular
concordance. For null-homotopic knots, a slightly weaker equivalence
relation is shown to admit a group structure.}



\maketitle

\section{Introduction}

Working in the setting of finite type invariants, Kirk and Livingston
defined families of type 1 invariants for knots in many 3--manifolds
and described indeterminacies that may arise when one attempts to
define relative invariants for homotopically essential pairs of knots
\cite{KL1, KL2}. The invariants of Kirk and Livingston are
extracted from the homology classes of double-point loops
corresponding to crossing changes in the knots during a homotopy. This
paper will extract invariants from the {\em homotopy} classes of such
double-point loops. The choices involved in identifying free homotopy
classes with elements in the fundamental group of the ambient manifold
require keeping careful track of conjugation actions and the resulting
linking ``numbers'' take values in a target orbit space which is in
general just a set (with a well-defined zero element) and depends in a
rather subtle way on the knots and manifolds being considered. Some of
the benefits justifying this somewhat unusual target include extending
the Kirk-Livingston invariants non-trivially to many more manifolds
(e.g.\ integral homology spheres) and showing they are concordance (not
just isotopy) invariants, as well as providing a clear geometric
characterization of the invariants in terms of Whitney disks and
demonstrating computability of the indeterminacies.

The combination of 3-- and 4--dimensional methods used in \cite{KL1} and \cite{KL2} highlights the
suggestive relation between the various crossing change diagrams of 3--dimensional finite type
theory and cross-sections of generic singularities of surfaces in 4--dimensions. The central idea
of this paper is that this relation can be further exploited by applying (a generalization of)
Wall's quadratic intersection form to the trace of a homotopy of knots. This approach provides a
connection between two advancing lines of research: the extensions of finite type theories to
arbitrary 3--manifolds (e.g.\ \cite{CM, GL, K, KL1, KL2} and many
others) and the developing theory of {\em Whitney towers}  which detect the failure of the Whitney
move in dimension 4 (e.g.\ \cite{COT, CST, S1, S3, ST, ST2}).

Knot theory in a non-simply connected manifold breaks naturally into the study of free homotopy
classes of knots (or classes of {\em 0-equivalent} knots in the language of finite type theory). We
assume that our manifolds are oriented and equipped with basepoints. For each element
$\gamma\in\pi_1M$ in the fundamental group of a 3--manifold $M$, let $\mathcal{K}_{\gamma}(M)$
denote the set of oriented knots in the free homotopy class determined by (the conjugacy class of)
$\gamma$ and let $\mathcal{C}_{\gamma}(M)$ denote the set $\mathcal{K}_{\gamma}(M)$ modulo
concordance (details in \ref{conc-singconc-subsec}). Thus, any two knots $k$ and $j$ in
$\mathcal{K}_{\gamma}$ co-bound an {\em immersed} annulus in $M\times I$ ($I$ the unit interval)
and our invariants will (1) provide obstructions to $k$ and $j$ representing the same element in
$\mathcal{C}_{\gamma}(M)$, that is, obstructions to $k$ and $j$ co-bounding an {\em embedded}
annulus in $M\times I$, and (2) provide geometric information towards improving a {\em singular}
concordance in the sense that all singularities can be paired by Whitney disks.

Although we will not work explicitly in the setting of finite type invariants, the reader familiar with the type 1
invariants of Kirk and Livingston will notice that they factor through the invariants defined here (with their
chosen cohomology class corresponding to a representation from $\pi_1M$ to a cyclic group).

A survey of our main results follows.
\paragraph{Null-homotopic knots} Our first observation
(details in Section~\ref{null-homotopic-sec}) is that, for any (oriented null-homotopic) knot
$k\in\mathcal{K}_1(M)$, Wall's self-intersection invariant $\mu$ (which counts signed double-point
loops) can be applied to the trace of any null-homotopy of $k$ in $M\times I$, yielding an {\em
algebraic self-linking number} $\mu(k)$ which takes values in a quotient $\widetilde{\Lambda}$ of
the free abelian group $\Lambda:=\Z[\pi_1M]$ generated by the elements of $\pi_1M$. The quotient is
by two relations: The inversion relation $g=g^{-1}$ for group elements corresponds to changing the
choice of orientation of a double-point loop, and the cusp relation $1=0$ corresponds to the fact
that local cusp homotopies create or eliminate double point loops having the trivial group element
$1\in\pi_1X$ (we are working with unframed knots). The concordance invariance of self-linking
numbers is a direct consequence of the 3--dimensional Sphere Theorem (see
Lemma~\ref{embedded-sphere-lemma}) and by clasp-doubling embedded loops it is not hard to construct
knots realizing all elements in $\widetilde{\Lambda}$; this is the content of the following theorem
which is proved in Section~\ref{null-homotopic-sec}.
\begin{thm}\label{absolute-self-link-thm}
The map $k\mapsto \mu(k)$ induces a well-defined map from $\mathcal{C}_1(M)$ onto $\widetilde{\Lambda}$.
\end{thm}
Due to the choice of basing of $k$, the image of $\mu(k)\in\widetilde{\Lambda}$ is only
well-defined up to conjugation by elements of $\pi_1M$; however, there are no further
indeterminacies and if $M$ is not simply connected then self-linking numbers detect many
null-homotopic knots which are not null-concordant, i.e., which do not bound an {\em embedded}
2--disk in $M\times I$. For example, if a knot $k_g\in M$ is constructed as the clasped-double of a
null homologous loop which represents any non-trivial element $g$ in the commutator subgroup of
$\pi_1M$, then $\mu(k_g)=\pm g\neq 0\in\widetilde{\Lambda}$ (Figure~\ref{clasp-knot-fig} in
\ref{onto-subsubsection}) showing that $k_g$ is not null-concordant and not 1-equivalent to the
unknot (see paragraph after Theorem~\ref{spherical-thm} below), facts not detected by the
homological invariants of \cite{KL1} and \cite{KL2}.

\paragraph{Essential knots}
In Section~\ref{rel-self-link-sec} we will define {\em relative} algebraic self-linking numbers for
pairs of {\em essential} knots, i.e.\ knots in $\mathcal{K}_{\gamma}(M)$ for $\gamma\neq 1$, by
applying a generalization of Wall's self-intersection invariant to immersed annuli in $M\times I$.
This generalized $\mu$--invariant for immersed annuli naturally takes values in an abelian group
$\tilde{\Lambda}_{\gamma}$ generated by the double cosets of the fundamental group of the ambient
manifold by the cyclic subgroup $\langle \gamma \rangle$ generated by the annulus; in our case,
$\gamma$ will be a knot longitude. If $k$ and $j$ are knots in $\mathcal{K}_{\gamma}(M)$ then the
{\em relative self-linking number $\mu_k(j)$ of $j$ with respect to $k$} takes values in an orbit
space of $\tilde{\Lambda}_{\gamma}$ under two group actions: There is a conjugation action by the
centralizer $\zeta(\gamma)$ which corresponds to basing choices and a more subtle action of a
certain {\em indeterminacy sub-group}
$$
\Phi(k)\leq {\tilde{\Lambda}}_{\gamma}\rtimes \zeta(\gamma)
$$
which accounts for the effect that pre-composing by a singular self-concordance of $k$ has on
$\mu_k(j)$.

\begin{figure}[ht!]
         \centerline{\includegraphics[scale=0.55]{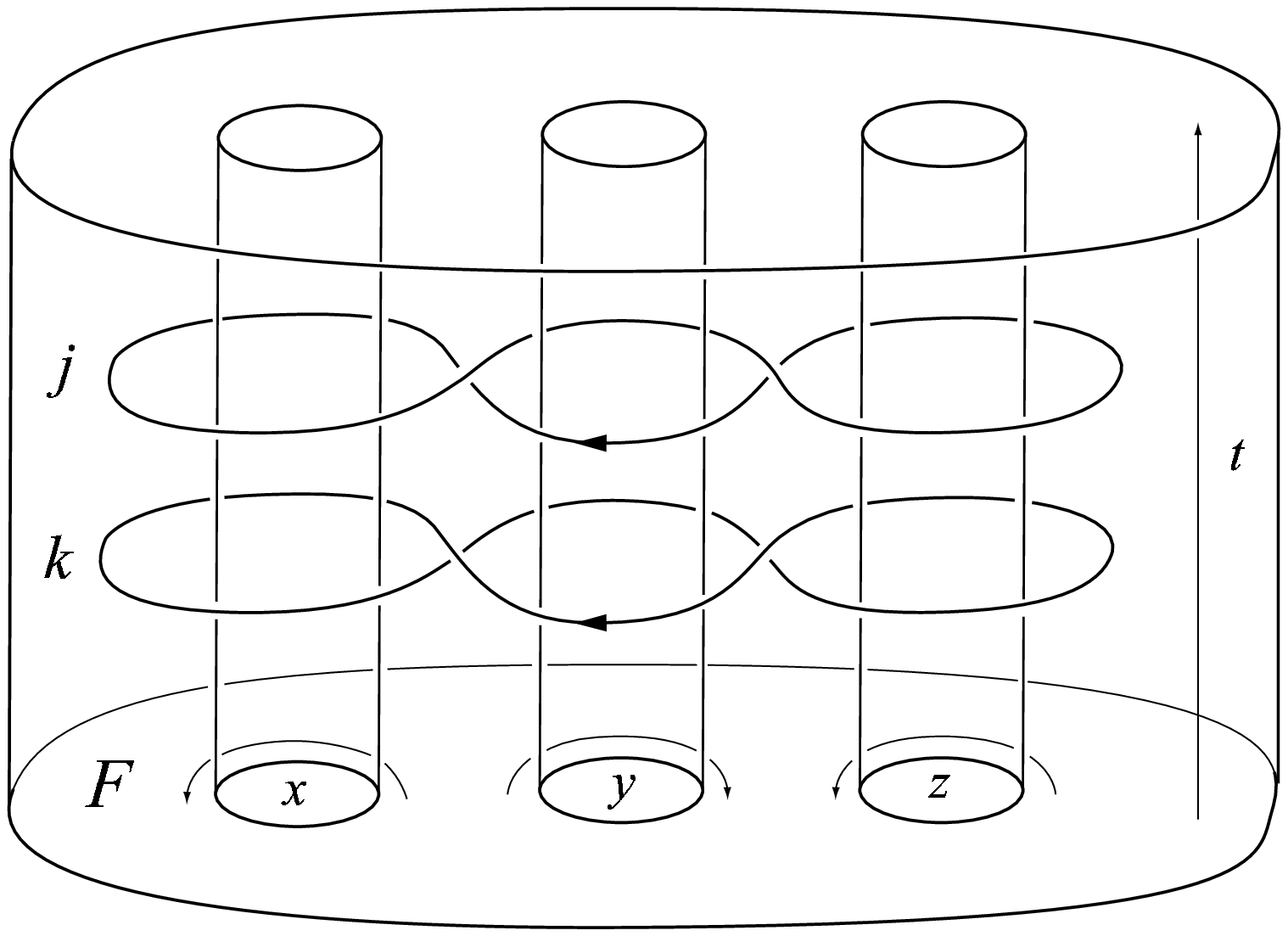}}
         \nocolon \caption{}
         \label{rel-example1-fig}

\end{figure}

The following example illustrates some basic properties of relative self-linking numbers. More
subtle properties will be exhibited in Section~\ref{example-sec}.
\paragraph{Example}
Illustrated in Figure~\ref{rel-example1-fig} are two knots $k$ and $j$ in $\mathcal{K}_{xyz}(M)$
where $M$ is the product $F\times S^1$ of a thrice punctured 2--disk $F$ with the circle and
$\pi_1M=\langle x,y,z \rangle\times\langle t \rangle$ is the cartesian product of the free group on
$x$, $y$ and $z$ (represented by loops in $F$ around the punctures) with the central cyclic group
generated by the element $t$ (represented by a circle factor of $M$). The figure shows $M$ cut open
along $F$ cross a point.

In this case, with respect to $k$, both the conjugacy and indeterminacy subgroup actions are
trivial (see Example~\ref{rel-example1-trivial}) so the relative self-linking number $\mu_k(j)$
takes values in the group $\widetilde{\Lambda}_{xyz}$. The two obvious crossing changes that change
$k$ into $j$ describe the trace $H$ of a homotopy with two singularities of differing sign. The
group elements associated to these singularities can be computed from the figure yielding (up to a
sign from global orientation conventions) $\mu_k(j)=\mu(H)=x-xy$ which is non-zero in
$\widetilde{\Lambda}_{xyz}$. That $k$ and $j$ are not isotopic can also, in this case, be detected
by the Kirk-Livingston invariants (using a cohomology class that is non-trivial on $y$); as a
result of the next theorem we see, furthermore, that $k$ and $j$ are not concordant.

Although the target for $\mu_k(j)$, denoted ${\tilde{\Lambda}}_{\gamma}/(\Phi(k),{\zeta(\gamma)})$,
does not in general have a group structure (Remark~\ref{no-group-remn}), it does have a
well-defined orbit of $0\in\tilde{\Lambda}_{\gamma}$ and so it makes sense to speak of $\mu_{k}(j)$
vanishing. The following theorem will be proved in Section~\ref{rel-self-link-sec}.
\begin{thm}\label{rel-self-link-thm}
For each $k\in\mathcal{K}_{\gamma}(M)$, the map $j\mapsto \mu_{k}(j)$ vanishes on knots concordant
to $k$ and induces a well-defined map
$$
\mathcal{C}_{\gamma}(M)\twoheadrightarrow{\tilde{\Lambda}}_{\gamma}/(\Phi(k),{\zeta(\gamma)})
$$
onto the target.
\end{thm}

Thus, the target reflects some of the structure of $\mathcal{C}_{\gamma}(M)$ and a variety of
properties which can occur are illustrated in the examples of Section~\ref{example-sec} (including
torsion (\ref{S^1-S^2})). Describing this target depends to a large extent on understanding the
indeterminacy sub-group $\Phi(k)$.  A satisfying theme that will emerge is that these
indeterminacies, which come from essential tori and 2--spheres in $M$, can be measured in terms of
intersections between lower-dimensional manifolds. A generating set for $\Phi(k)$ is given in
Proposition~\ref{Phi-prop}.

The indeterminacies due to 2--spheres are computed in terms of an intersection pairing
$\widetilde{\lambda}(\sigma,k)$ between $\pi_2(M)$ and $\mathcal{K}_{\gamma}(M)$
(\ref{sphere-ints}). This allows our knot invariants to be defined in reducible 3--manifolds, which
is desirable even if one is ultimately interested in the manifold structure; for instance,
surgering a knot in a reducible manifold can yield an irreducible one.

The indeterminacies due to tori correspond to singular self-concordances which project to essential
immersed tori in $M$ and only exist if the centralizer $\zeta(\gamma)$ of $\gamma$ is non-cyclic.
In this case, the computation of $\Phi(k)$ can be reduced to a computation in a Seifert fibered
submanifold that carries $\zeta(\gamma)$ using the 3--manifold structure theorems of Jaco-Shalen
\cite{JS1, JS2}, Johannson \cite{J}, Gabai \cite{G} and Casson-Jungreis \cite{CJ}. This
follows the approach used by Kalfagianni \cite{K} to show the existence of finite type knot
invariants in many irreducible 3--manifolds. These indeterminacies correspond to intersections
between knots and tori in $M$ and can in general be computed in terms of an intersection pairing
between curves in an orbit surface of a Seifert fibered sub-manifold of $M$ (\cite{S2}). In this
paper we will mostly be concerned with cases where these indeterminacies are trivial.

\paragraph{Spherical knots} It turns out
that there exist many knots for which the indeterminacy subgroup $\Phi(k)$ is completely determined
by the pairing $\widetilde{\lambda}(\sigma,k)$ between the knots and 2--spheres; these knots,
called {\em spherical knots} (\ref{spherical-subsec}), play a preferred role (analogous to the
unknot in the trivial homotopy class) in the geometric interpretation of linking invariants:
 \begin{thm}\label{spherical-thm}

 If $k\in\mathcal{K}_{\gamma}(M)$ is spherical, then for any knots $j$ and $j'$
 in $\mathcal{K}_{\gamma}(M)$ the
 following are equivalent:
 \begin{enumerate}
  \item $\mu_k(j)=\mu_k(j')$.
  \item  There
exists a singular concordance between $j$ and $j'$ such that all singularities are paired by Whitney disks.
  \end{enumerate}

\end{thm}
As illustrated later in Example~\ref{not-spherical-knot-example}, there are examples where
$\mu_k(j)=\mu_k(j')$ for a {\em non}-spherical knot $k$ and property (ii) of
Theorem~\ref{spherical-thm} does {\em not} hold (as detected by $\mu_{k_0}(j)\neq\mu_{k_0}(j')$ for
some spherical $k_0$).

\paragraph{$W\!$--equivalence}\label{W-equiv-intro}
Theorem~\ref{spherical-thm} shows that, with respect to spherical knots, relative self-linking
numbers inherit the geometric characterization of Wall's self-intersection invariant in terms of
Whitney disks (see Proposition~\ref{wall-int-prop}). The equivalence relation on ${K}_{\gamma}(M)$
defined by property (ii) of Theorem~\ref{spherical-thm} is called {\em $W\!$--equivalence} and
represents a ``first order'' improvement over an arbitrary singular concordance (a 0-equivalence)
in the following sense: Conant and \mbox{Teichner} have characterized finite type $n$--equivalence
(in the sense of Habiro and Goussarov) of knots in 3--manifolds in terms of the notion of
3--dimensional {\em capped grope cobordism} \cite{CT1, CT2}. In their language,
1--equivalence corresponds to capped surface cobordism. Such a capped surface bordism in $M$ can be
pushed into $M\times I$ and surgered to a $W\!$--equivalence. (See \cite{S1} for details including
higher orders.) Thus, relative self-linking numbers give obstructions to 1--equivalence of knots.

\paragraph{Example} In the example illustrated in
Figure~\ref{rel-example1-fig} above, the knot $k$  is spherical (see examples
\ref{rel-example1-trivial} and \ref{rel-example1-spherical}), hence $k$ and $j$ are not
$W\!$--equivalent, by taking $j'=k$ in Theorem~\ref{spherical-thm}, since
$\mu_k(j)-\mu_k(k)=\mu_k(j)$ is non-zero.\medskip

We are assured of the existence of spherical knots in {\em all} homotopy classes by restricting to a
large class $\mathcal{M}$ of 3--manifolds which do not contain circle bundles over non-orientable
surfaces whose total spaces are orientable and do not contain certain Seifert fibered spaces
containing non-vertical tori as sub-manifolds (see Section~\ref{M-sec}).

\begin{thm}\label{M-thm}
For any $M\in\mathcal{M}$ and any $\gamma\in\pi_1M$, there exists a spherical knot $k_0$ in
$\mathcal{K}_{\gamma}(M)$.
\end{thm}
Thus, relative self-linking numbers characterize $W\!$--equivalence in $\mathcal{M}$.

\paragraph{2--component links}

The discussion so far can be applied similarly to define (relative) concordance invariants for
2--component links in terms of (a generalization (\ref{annuli-ints}) of) Wall's intersection
pairing $\lambda$ (\ref{wall-form}) with analogous results. This is sketched in
Section~\ref{link-sec}. One result worth mentioning here is that the (absolute) {\em algebraic
linking number} (\ref{absolute-link-defi}) for a 2--component link of null-homotopic knots plays a
role in defining a group structure on a quotient of $K_1(M)$ (see \ref{group-cor}).

\paragraph{Conventions}

For the most part, standard 3-- and 4--dimensional techniques and terminology are used throughout.
Irreducible 3--manifolds will be allowed to have spherical boundary components. The closed unit
interval $[0,1]$ will be denoted by $I$ and occasionally be reparametrized implicitly. We work in
the smooth oriented category with specific orientations usually suppressed.

\paragraph{Acknowledgments} 
I am happy to thank Paul Kirk and Peter Teichner for helpful
conversations, and my former advisor Rob Kirby for his guidance and
support. Thanks also to the referee whose careful reading and
thoughtful comments have significantly contributed to improving the
exposition. This work was supported in part by an NSF Postdoctoral
Fellowship and the Max-Planck-Institut f\"{u}r Mathematik.


\section{Preliminaries}\label{prelim-sec}
This section briefly reviews the 4--dimensional version of Wall's intersection and
self-intersection invariants, $\lambda$ and $\mu$, as well as the notions of concordance and
singular concordance and also serves to fix notation. See also \cite{FQ} for more details.

\subsection{Wall's intersection invariants}\label{wall-form}
Let $D$ and $E$ be {\em properly immersed} 2--spheres or 2--disks (rel $\partial$) in a 4--manifold
$X$, that is, boundary is embedded in boundary and interior immersed in interior. After a small
perturbation (rel $\partial$), $D$ and $E$ can be assumed to be in general position, so that $D$
meets $E$ in a finite set of transverse {\em intersection points} and each of $D$ and $E$ have
finitely many transverse {\em self-intersection points} or {\em double point} singularities.
Neighborhoods of an intersection point $p$ in $D$ and $E$ are called {\em sheets} of $D$ and $E$
(at $p$). Fix {\em whiskers} for each of $D$ and $E$; that is, choose an arc in $X$ connecting the
basepoint of $X$ to a basepoint on $D$ and likewise for $E$.

\subsubsection{Intersection numbers}
Each point $p\in D\cap E$ determines an element $g_p \in \pi_1X$ from the following loop: First go
along the whisker on $D$ from the basepoint of $X$ to the basepoint of $D$, then along $D$
(avoiding all double points) to $p$, then along $E$ (avoiding all double points) to the basepoint
on $E$ and then back along $E$'s whisker to the basepoint of $X$. Since $D$ and $E$ are simply
connected, $g_p$ does not depend on how the loop runs between $p$ and the basepoints on $D$ and
$E$. By summing (with appropriate signs) over all intersection points we get an intersection
``number'' in $\Lambda:=\Z[\pi_1X]$, the free abelian group generated by the elements of $\pi_1X$:

\begin{defi}

The {\em intersection number} $\lambda(D,E)$ of $D$ and $E$ is defined by
$$
\lambda(D,E):=\sum (\sign  p)\cdot g_p \in \Lambda
$$
where the sum is over all intersection points $p\in D\cap E$ and
$\sign  p$ equals $+1$ (resp.\ $-1$) if the orientation of $X$ at
$p$ agrees (resp.\ disagrees) with the orientation determined by
the sheets of $D$ and $E$ at $p$.
\end{defi}

Note that the $g_p$ are all computed using the fixed whiskers on
$D$ and $E$. Changing the whisker on $D$ (resp.\ $E$) changes
$\lambda(D,E)$ by left (resp.\ right) multiplication by an element
of $\pi_1X$.

\subsubsection{Self-intersection numbers}

\begin{figure}[ht!]
         \centerline{\includegraphics[scale=.55]{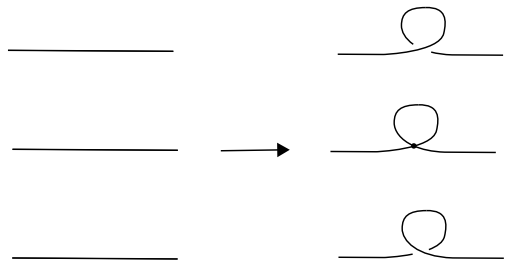}}
         \caption{Before and after a local cusp homotopy of a surface in a 4--manifold}
         \label{cusp-fig}

\end{figure}

For each double point $p$ of $D$, define $g_p\in \pi_1X$ from the following loop: First go along
the whisker on $D$ from the basepoint of $X$ to the basepoint of $D$, then along $D$ (avoiding all
double points) to $p$, then change sheets at $p$ and go back along $D$ (avoiding all double points)
to the basepoint on $D$ and then back along $D$'s whisker to the basepoint of $X$. Note that (for a
fixed whisker) $g_p$ depends only on the choice of first sheet at $p$ and changing the order of
sheet-change for the loop at $p$ changes $g_p$ to ${g_p}^{-1}$. Note also that a local cusp
homotopy (Figure~\ref{cusp-fig}) creates a double point $p$ with $g_p$ equal to the trivial element
$1\in \pi_1X$.

\begin{defi}\label{mu-defi}
The self-intersection number $\mu(D)$ is defined by
$$
\mu(D):=\sum (\sign  p)\cdot g_p\in \tilde{\Lambda}
$$
where the sum is over all double points $p$ of $D$ and $\sign  p =\pm1$ is determined by comparing the orientation
of $X$ at $p$ with the orientation given by the two sheets of $D$ at $p$.
\end{defi}
Here $\tilde{\Lambda}$ is as described in the introduction:
$$
\tilde{\Lambda}:=\frac{\Lambda}{\{g-g^{-1}\}\oplus\Z[1]}
$$
where $g$ ranges over $\pi_1X$ and $\Z[1]$ is generated by the trivial element $1\in \pi_1X$. (This
is a quotient as an abelian group, not as a ring.) Note that changing the whisker on $D$ changes
$\mu(D)$ by conjugation by an element of $\pi_1X$.

\paragraph{Remark}This definition of self-intersection number is sometimes referred to as the {\em
reduced} self-intersection number; omitting the quotient by $\Z[1]$ yields an unreduced version
which is only invariant under {\em regular} homotopy. Using the unreduced version would lead to an
invariant of {\em framed} knots.

\subsubsection{Homotopy invariance and the geometry of Wall's
invariant}

The homotopy invariance of intersection numbers can be shown by an
analysis of the singularities of homotopies of surfaces
(\cite{FQ}):
\begin{prop}\label{wall-int-prop}
In the above setting, $\lambda(D,E)\in \Lambda$ and $\mu(D)\in \tilde{\Lambda}$ depend only on the homotopy
classes (rel $\partial$) of $D$ and $E$. The intersection number $\lambda(D,E)$ and self-intersection numbers
$\mu(D)$ and $\mu(E)$ all vanish if and only if $D$ and $E$ are homotopic (rel $\partial$) to disjoint maps with
all self-intersections paired by Whitney disks.\endproof
\end{prop}
The disjointness property follows from standard manipulations of Whitney disks. (See 1.4 of
\cite{FQ} for an introduction to Whitney disks in 4--manifolds.) If the Whitney disks in the
conclusion of Prop~\ref{wall-int-prop} happened to be disjointly embedded, correctly framed and
with interiors disjoint from $D$ and $E$, then it would follow that $D$ and $E$ are homotopic (rel
$\partial$) to disjoint embeddings.

Note that in a simply connected 4--manifold $\lambda(D,E)$ reduces to the usual homological
intersection $D\cdot E \in \Z$ which counts the signed intersection points and that
$\tilde{\Lambda}\cong {0}$ so that $\mu$ always vanishes.

\subsection{Concordance and singular concordance}\label{conc-singconc-subsec}

\begin{defi}\label{conc-defi}
Two knots $k$ and $j$ in a 3--manifold $M$ are {\em concordant} if there exists a properly embedded
annulus $(C,\partial C) \hookrightarrow (M\times I,M\times \partial I)$ bounded by $k \subset M
\times \{0\}$ and $j \subset M\times \{1\}$. The oriented boundary of the annulus $C$ is required
to induce the difference of the orientations of the knots: $\partial C= k-j$. Such a $C$ is called
a {\em concordance} between $k$ and $j$. Two links (collections of disjoint knots) $l$ and $l'$ in
$M$ are {\em concordant} if their components can be joined by a collection of properly embedded
pairwise disjoint annuli in $M\times I$ where each annulus is a concordance between a component of
$l$ and a component of $l'$. The union of the annuli is a {\em concordance} between $l$ and $l'$.
\end{defi}

Concordance is clearly an equivalence relation on knots or links with a fixed number of components
and isotopy is a special kind of level-preserving concordance. By allowing the annuli to be
immersed we get the (much) weaker equivalence relation of {\em singular concordance}:

\begin{defi}\label{sing-conc-defi}
Two knots $k$ and $j$ in $M$ are {\em singularly concordant} if
there exists a properly immersed annulus $(C,\partial C)
\looparrowright (M\times I,M\times \partial I)$ bounded by $k
\subset M \times \{0\}$ and $j \subset M\times \{1\}$. Such a $C$
is called a {\em singular concordance} between $k$ and $j$ (or
from $k$ to $j$). Two links $l$ and $l'$ in $M$ are {\em
singularly concordant} if their components can be joined by a
collection of properly immersed annuli (not necessarily disjoint)
in $M\times I$ where each annulus is a singular concordance
between a component of $l$ and a component of $l'$.
\end{defi}

The trace of a homotopy of knots is a (level-preserving) singular concordance and, by \cite{Go} and
\cite{Gi}, the relations of singular concordance and homotopy are in fact equivalent for knots and
links in 3--manifolds. It will be convenient to make statements in the a priori more general
language of singular concordance, however explicit constructions will usually be described by
homotopies.

\subsection{Whiskers for knots}\label{knot-whiskers}
Recalling the notation in the introduction, we will use the fundamental group to index singular
concordance classes of knots:
\begin{defi}\label{conj-class-defi}
For each element $\gamma$ in the fundamental group of a 3--manifold $M$, let
$\mathcal{K}_{\gamma}(M)$ denote the set of oriented knots (up to isotopy) in the free homotopy
class determined by (the conjugacy class of) $\gamma$ and let $\mathcal{C}_{\gamma}(M)$ denote the
set $\mathcal{K}_{\gamma}(M)$ modulo concordance.
\end{defi}
This means that whenever we connect a knot $k\in\mathcal{K}_{\gamma}(M)$ by a whisker to the basepoint of $M$ we
require that this basing satisfies $[k]=\gamma\in\pi_1M$.

\section{Null-homotopic knots}\label{null-homotopic-sec}
This section contains the precise definition of the algebraic self-linking number $\mu(k)$ and the
proof of Theorem~\ref{absolute-self-link-thm}. The arguments here also apply to algebraic linking
numbers of 2--component links of null-homotopic knots as will be described in
Section~\ref{link-sec}.

\subsection{Algebraic self-linking numbers}
Recall that our knots and manifolds are assumed oriented.
\begin{defi}\label{absolute-self-link-defi}
For $k\in\mathcal{K}_{1}(M)$ define the {\em algebraic self-linking number} $\mu(k)$ by
$$
\mu(k):=\mu(D)\in\widetilde{\Lambda}
$$
where $D$ is any properly immersed 2--disk in $M\times I$ bounded by $k\subset M\times\{0\}$ and
$\mu(D)$ is Wall's self-intersection number as defined in \ref{mu-defi}. Here $\pi_1(M\times I)$ is
identified with $\pi_1M$ via projection onto $M\times\{0\}$. The immersed 2--disk $D$ is a {\em
singular null-concordance} of $k$ and is oriented by the orientation of $k$ (via some fixed
convention).
\end{defi}

We may refer to $\mu(k)$ simply as a ``self-linking number'' for sake of brevity, since in the
present context omitting ``algebraic'' should not cause confusion with the well known
integer-valued self-linking number of a framed knot.

\subsection{Proof of Theorem~\ref{absolute-self-link-thm}}

\begin{figure}[ht!]
         \centerline{\includegraphics[scale=.45]{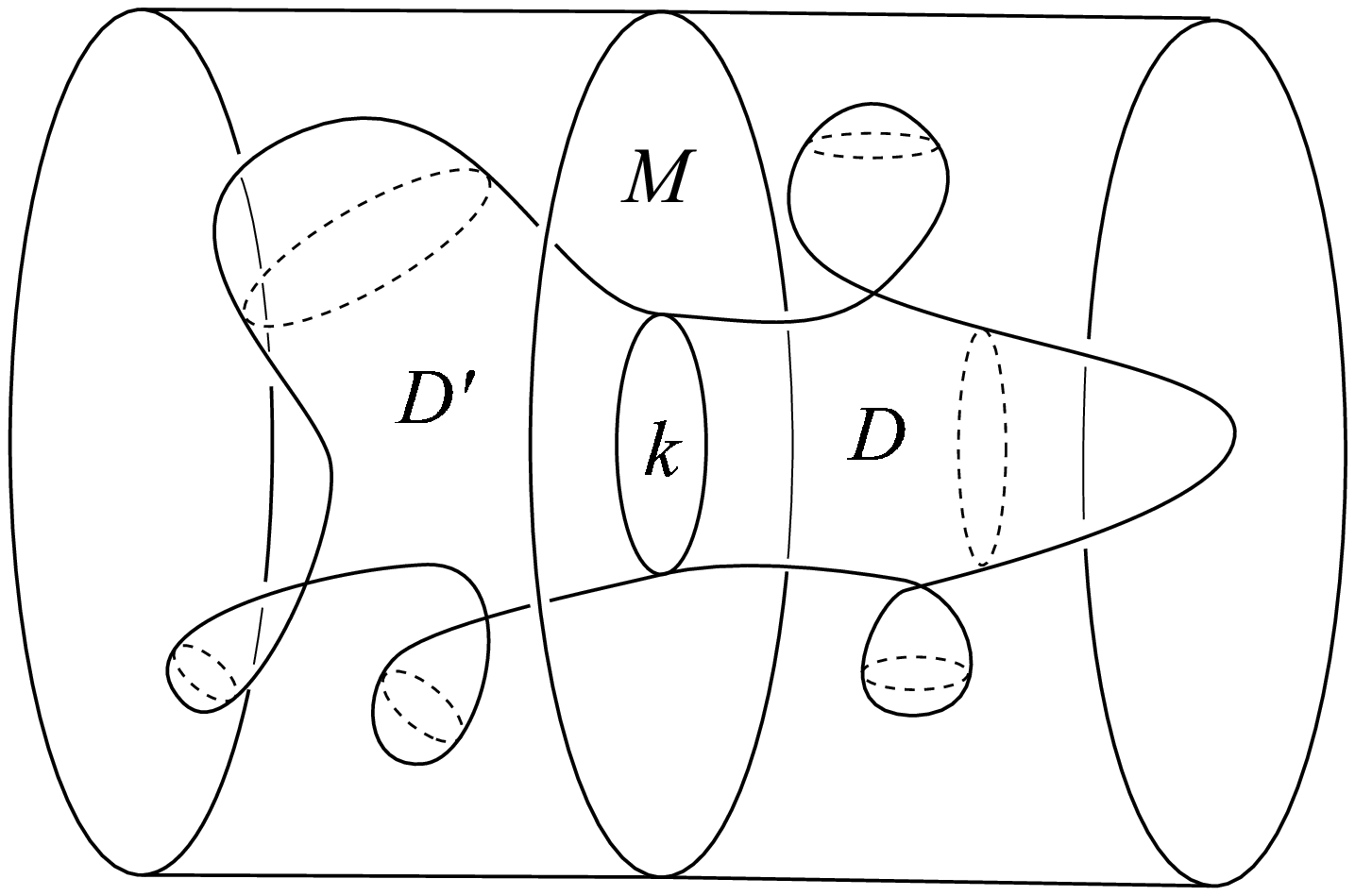}}
         \nocolon \caption{}
         \label{2disks-fig}

\end{figure}

\paragraph{Independence of choice of $D$} 
We show first that $\mu(k)$ does not depend on the choice of singular
null-concordance. Let $D$ and $D'$ be two singular null-concordances
of $k$. Then the union $S$ of $D$ and $D'$ along $k$ (in two copies of
$M\times I$ identified along $M$) determines an element of
$\pi_2(M\times I)\cong \pi_2(M)$ (see Figure~\ref{2disks-fig}).

\begin{lem}\label{embedded-sphere-lemma}
Any $n$ elements of $\pi_2(M\times I)$ are represented by $n$ embedded pairwise disjoint
2--spheres.
\end{lem}

\begin{proof}
A well known consequence of the 3--dimensional Sphere Theorem is that $\pi_2(M)$ is generated as a
module over $\pi_1(M)$ by disjoint embeddings (the 2--spheres that decompose $M$ into prime
factors, together with any spherical boundary components and cross-sections of any $S^2 \times S^1$
factors, see Proposition 3.12 of \cite{Ha}). Tubing these generators together in $M\times I$ does
not create any new intersections, so $\pi_2(M\times I)$ is spanned by disjoint embeddings.
\end{proof}

Lemma~\ref{embedded-sphere-lemma} implies in particular that
Wall's intersection form vanishes on $\pi_2(M\times I)$ and so
$$
\mu(S)=0=\mu(D)-\mu(D').
$$
This shows that $\mu(k)$ does not depend on the choice of bounding disk.

\paragraph{Concordance invariance}
If $C$ is a concordance from $k'$ to $k$ then, up to conjugation, $\mu(k')=\mu(C\cup D)=\mu(D)=\mu(k)$ since $C$
has no singularities.

\begin{figure}[ht!]
         \centerline{\includegraphics[scale=.45]{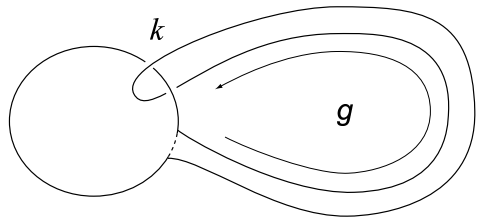}}
         \caption{Creating a null-homotopic knot $k$ with $\mu(k)=\pm g$}
         \label{clasp-knot-fig}

\end{figure}

\paragraph{The map $k\mapsto\mu(k)$ is
onto}\label{onto-subsubsection} To construct a null-homotopic knot $k_g \subset M$ with
$\mu(k_g)=\pm g$ for any element $g\in \pi_1(M)$, push an arc of a small circle around a loop
representing $g$ and create a $\pm$-clasp with the circle (Figure~\ref{clasp-knot-fig}). By
iterating this procedure (or band summing together such clasps), one can realize anything in
$\Lambda$ as $\mu(k)$ for some $k$. \endproof

\begin{rem}\label{null-remn}
Note that Lemma~\ref{embedded-sphere-lemma} implies that {\em any} homotopy invariant of 2--spheres
in a 4--manifold which vanishes on disjoint embeddings (and restricts to an invariant of immersed
2--disks rel boundary) will give a concordance invariant of null-homotopic knots in a 3--manifold.
In particular, the invariant $\tau$ of \cite{ST} can be used to define second order knot invariants
whenever the first order invariants vanish (see \cite{S}).
\end{rem}


\section{Relative self-linking
numbers}\label{rel-self-link-sec}

This section introduces the definitions, conventions and notation necessary to define relative
self-linking numbers $\mu_k(j)$ and gives the proof of Theorem~\ref{rel-self-link-thm}. Since
basepoint issues will complicate the procedure, some motivation is in order. The goal is to
construct a relative knot invariant that is as free from basepoint dependence as possible {\em and}
has a clear geometric characterization. The first step (\ref{annuli-self-ints}) is to generalize
Wall's $\mu$ homotopy invariant to immersed annuli in such a way that $\mu(A)$ vanishes if and only
if the singularities of the immersed annulus $A$ can be paired up by Whitney disks (after a
homotopy rel $\partial$). This is where the double coset space (by the cyclic subgroup $\langle
\gamma \rangle$) enters the picture. Here also, the invariant in $\tilde{\Lambda}_{\gamma}$ is
well-defined only up to a conjugation action by the centralizer subgroup $\zeta(\gamma)$,
corresponding to the choice of whisker connecting the annulus to the basepoint of the ambient
manifold.

Next, we want to apply $\mu$ to an immersed annulus joining two knots in $M\times I$. (Such an
immersed annulus will often be denoted $H$ since in practice it will usually be the trace of a
homotopy between knots.) This requires fixing conventions (\ref{conc-conventions}) and more
importantly examining the behavior of $\mu$ under composition of singular concordances, in
particular singular {\em self}-concordances of a knot which may contribute indeterminacies to the
final invariant. A subtle effect that arises is {\em another} conjugation action by $\zeta(\gamma)$
which comes from a loop traced out by a basepoint on the knot during a homotopy, together with the
whiskers on the knots at either end of the homotopy; such a loop is called a {\em latitude}
(\ref{conc-latitudes}, \ref{self-conc-latitudes}) and this conjugation action prevents $\mu$ from
being additive (\ref{conc-ints}, \ref{composing-concs}). The effect of this non-additivity is
captured by the action of an {\em indeterminacy subgroup} $\Phi(k)$ which depends in general on the
``base'' knot $k$ (\ref{Phi(k)-subsec}). This action is compatible with the basepoint conjugation
action (\ref{conj-action-subsec}) and the precise definition of $\mu_k(j)$ is given in
\ref{rel-self-link-subsec}.

Having successfully waded through these conventions and notations, the proof of
Theorem~\ref{rel-self-link-thm} then more or less falls out of the definitions. In digesting the
above conjugation actions, notations, etc., it may be helpful for the reader to draw schematic
pictures along the lines of Figure~\ref{H-conc-fig} (in \ref{rel-self-link-thm-proof}) in order to
see the effect of compositions, whisker changes, etc. on the double-point loops counted by $\mu$.

\subsection{Self-intersection numbers for
annuli}\label{annuli-self-ints} Let $A:(S^{1}\times I, S^{1} \times \{0,1\})\looparrowright
(X,\partial X)$ be a properly immersed annulus in a 4--manifold $X$. Choosing a whisker for $A$
identifies the image of $\pi_1(S^1 \times I)\cong \Z$ with a cyclic subgroup $\langle \gamma
\rangle$ of $\pi_1X$ generated by the image of a generating circle of the annulus. For each double
point $p$ of $A$, associate an element $g_p\in \pi_1X$ from a loop in $A$ that changes sheets at
$p$ (and avoids all other double points) together with the whisker on $A$. Note that $g_p$ is only
well-defined up to left and right multiplication by powers of $\gamma$, corresponding to the loop
in $A$ wandering around the circle direction before and after changing sheets at $p$. Denote by
$\Lambda_{\gamma}$ the free abelian group additively generated by the double cosets of $\pi_1X$ by
$\langle \gamma \rangle$, that is,
$$
 \Lambda_{\gamma}:=\Z[{\langle \gamma \rangle}\backslash
\pi_1X/{\langle \gamma \rangle}] =\frac{\Z[\pi_1X]}{\{g-{\gamma}^n
g {\gamma}^m\}}
$$
where $n$ and $m$ range over the integers.

To account for indeterminacies due to the choices of orientation of the double-point loops and from
the introduction of trivial group elements via cusp homotopies, we take a further quotient and
define
$$
\tilde{\Lambda}_{\gamma}:=\frac{\Lambda_{\gamma}}{\{g-g^{-1}\}\oplus
\Z[1]} =\frac{\Z[\pi_1X]}{\{g-{\gamma}^n  g^{\pm 1}
{\gamma}^m\}\oplus \Z[1]}.
$$

\begin{defi}\label{annulus-int-defi}

Let $A$ be a properly immersed annulus in a 4--manifold $X$ and $\langle \gamma \rangle$ be the
image of the induced map on fundamental groups. Then the {\em self-intersection number} $\mu(A)$ of
$A$ is defined by
$$
\mu(A):=\sum (\sign  p)\cdot g_p \in \tilde{\Lambda}_{\gamma}
$$
where the sum is over all double points $p$ of $A$ and $\sign  p$ comes from the orientations of $X$ and the
sheets of $A$ at $p$ as usual.

\end{defi}
Changing the whisker for $A$ in a manner that preserves the image
$\langle \gamma \rangle$ of the fundamental group has the effect
of conjugating $\mu(A)$ by an element of the centralizer subgroup
$\zeta(\gamma):=\{g\in \pi_1X\mid g \gamma=\gamma g \} \leq
\pi_1X$ (see \ref{conc-ints} below). Thus, the orbit of $\mu(A)$
in $\tilde{\Lambda}_{\gamma}$ under conjugation by $\zeta(\gamma)$
is invariant.

In general, changing the whisker on $A$ changes both $\mu(A)$ and
the target space by an isomorphism:
$z\in\tilde{\Lambda}_{\gamma}\mapsto \alpha
z\alpha^{-1}\in\tilde{\Lambda}_{\alpha\gamma\alpha^{-1}}$.

The same arguments as in the proof of
Proposition~\ref{wall-int-prop} give:

\begin{prop}\label{annuli-selfint-prop}

The above defined self-intersection number $\mu(A)$ is invariant under homotopy (rel $\partial$) and vanishes if
and only if all singularities of $A$ can be paired by Whitney disks (perhaps after some local cusp
homotopies).\endproof
\end{prop}

\subsection{Conventions for singular concordances}\label{conc-conventions}

Let $H:(S^1 \times I,S^1 \times \{0,1\}) \rightarrow (M \times I,M \times \{0,1\})$ be a singular concordance
between any knots $k$ and $j$ in $\mathcal{K}_{\gamma}(M)$. Orientations of $H$ and $M\times I$ are determined by
the orientations on $k\subset M\times\{0\}$ and $M$ together with the orientation of $I$ so that $\partial H=k-j$.
Identify $\pi_1(M\times I)$ with $\pi_1M$ via projection onto $M\times \{0\}$. Take a whisker for $k$ as the
whisker for $H$. Since by convention (\ref{knot-whiskers}) we only take whiskers for $k$ so that $[k]=\gamma$, the
image of the fundamental group of the annulus equals $\langle \gamma\rangle\leq\pi_1M$.

Unless otherwise specified, these conventions will be assumed for
all singular concordances.

\subsection{Latitudes of singular concordances}\label{conc-latitudes}
Let $H$ be a singular concordance from $k$ to $j$ with whiskers chosen so that
$[k]=[j]=\gamma\in\pi_1M$. A {\em latitude} of $H$ is any arc that goes from the basepoint of $M$
in $M\times \{0\}$ along the whisker $w$ of $k$, then along $H$, then along the whisker $w'$ of $j$
to the basepoint of $M$ in $M\times \{1\}$. The projection of a latitude of $H$ to $M\times \{0\}$
is a loop that determines an element $\phi$ in the centralizer $\zeta(\gamma)$ of $\gamma$ in
$\pi_1M$. This element $\phi$ is well-defined up to multiplication by powers of $\gamma$ and hence
determines a well-defined element in the double coset space
$\langle\gamma\rangle\backslash\pi_1M/\langle\gamma\rangle$.

In this setting we will speak of ``{\em the} element $\phi\in\zeta(\gamma)$ determined by a latitude'' since such
an element represents a well-defined double coset.

The notation $H_{\phi}$ will sometimes be used to indicate that a
singular concordance has a latitude whose projection represents
$\phi$.

\subsection{Latitudes of singular self-concordances}\label{self-conc-latitudes}
A singular {\em self}-concordance $K$ of a knot $k$ with whisker
$w$ has a preferred latitude which uses $w$ at both ends of $K$. A
latitude of a singular self-concordance will always be assumed to
be such a preferred latitude and when it is necessary to carefully
specify the whisker being used we will use the notation
$$
\phi=\lat[K](w)
$$
for the element $\phi\in\zeta(\gamma)$ determined by $K$ and $w$.

\subsection{Self-intersection numbers of singular
concordances}\label{conc-ints} The self-intersection number $\mu(H)$ of a singular concordance $H$
of $k\in\mathcal{K}_{\gamma}(M)$ takes values in $\tilde{\Lambda}_{\gamma}$ and when necessary we
use the notation $\mu(H)(w)$ to indicate that $\mu(H)$ is computed using a specific whisker $w$ for
$k$. By our convention that $[k]=\gamma\in\pi_1M$, any two choices of whisker for $k$ differ by an
element $\alpha\in\zeta(\gamma)$, that is, $\alpha=[w'-w]$ is determined by the oriented loop
$w'-w$ which goes to $k$ along $w'$ and back to the basepoint of $M$ along $w$. The effect of a
whisker change on $\mu(H)$ is conjugation by $\alpha$:
$$
\mu(H)(w')=\alpha(\mu(H)(w))\alpha^{-1}.
$$

\subsection{Additivity and composition}\label{composing-concs}

If $H$ is a singular concordance from $k$ to $k'$ and $H'$ is a
singular concordance from $k'$ to $j$ then we write $H+H'$ for the
singular concordance from $k$ to $j$ that is the composition of
$H$ followed by $H'$ (with $I$ reparametrized appropriately). If
$w$ and $w'$ are the whiskers for $k$ and $k'$ we have
$$
\mu(H+H')=\mu(H)+\phi\mu(H')\phi^{-1}
$$
where $\phi$ is determined by a latitude of $H$ using the whiskers
$w$ and $w'$ (recall our convention that basepoints are taken in
$M\times \{0\}$). We will be particularly concerned with the case
where $k=k'$ and $w=w'$ so that $H$ is a singular self-concordance
and $\phi=\lat[H](w)$.

Note that if a latitude of $H'$ determines $\phi'$ then $H+H'$ has
a latitude determining $\phi\phi'$.

If $H$ is a singular concordance from $k$ to $j$ then we write $-H$ for the ``inverse'' singular
concordance from $j$ to $k$ that inverts the $I$ parameter of $S^1 \times I$ (which has the effect
of changing all the signs of the singularities of $H$). For instance, if $H$ is the trace of a
homotopy of a knot then $-H$ is the trace of the homotopy ``run backwards.'' If a latitude of $H$
determines $\phi$, then $-H$ has a latitude determining $\phi^{-1}$ and
$$
\mu(-H)=-\phi^{-1}\mu(H)\phi.
$$

\subsection{The indeterminacy subgroup $\Phi(k)$}\label{Phi(k)-subsec}

\begin{defi}\label{Phi-defi}

For any knot $k\in\mathcal{K}_{\gamma}(M)$, define
$$
\Phi(k)\leq {\tilde{\Lambda}}_{\gamma}\rtimes \zeta(\gamma)
$$
to be the subgroup of the semi-direct product of
${\tilde{\Lambda}}_{\gamma}$ and $\zeta(\gamma)$ (with respect to
the conjugation action of $\zeta(\gamma)$ on
${\tilde{\Lambda}}_{\gamma}$) generated by the elements
$$
(\mu(K)(w),\lat[K](w))\in{\tilde{\Lambda}}_{\gamma}\rtimes
\zeta(\gamma)
$$
as $K$ ranges over all singular self-concordances of $k$ in $M$
and $w$ ranges over all whiskers identifying
$[k]=\gamma\in\pi_1M$.

\end{defi}

\subsubsection{The action of $\Phi(k)$ on
${\tilde{\Lambda}}_{\gamma}$}

The group ${\tilde{\Lambda}}_{\gamma}\rtimes \zeta(\gamma)$ acts
on ${\tilde{\Lambda}}_{\gamma}$ by
$$
(z,\phi): y \mapsto z+\phi y \phi^{-1}
$$
for $(z,\phi)\in{\tilde{\Lambda}}_{\gamma}\rtimes \zeta(\gamma)$
and $y\in{\tilde{\Lambda}}_{\gamma}$. Denote the equivalence
classes under the restriction of this action to the subgroup
$\Phi(k)$ by
$$
{\tilde{\Lambda}}_{\gamma}/\Phi(k).
$$

The above action corresponds to the effect that pre-composing by a
singular self-concordance $K$ has on the self-intersection number
of a singular concordance $H$ (\ref{composing-concs}) {\em
provided that $\mu(K)$ and $\mu(H)$ are both computed using the
same whisker}.

\subsubsection{The action of $\zeta(\gamma)$ on
${\tilde{\Lambda}}_{\gamma}/\Phi(k)$}\label{conj-action-subsec}

The effect that changing the whisker for $k$ has on a generator of
$\Phi(k)$ is described by the diagonal conjugation action of
$\zeta(\gamma)$ on ${\tilde{\Lambda}}_{\gamma}\rtimes
\zeta(\gamma)$:
$$
\alpha:(z,\phi)\mapsto (\alpha z \alpha^{-1}, \alpha
\phi\alpha^{-1})
$$
for all $\alpha\in\zeta(\gamma)$.

We have a well-defined conjugation action of $\zeta(\gamma)$ on
${\tilde{\Lambda}}_{\gamma}/\Phi(k)$ since if
$$
x= z+\phi y \phi^{-1} \in{\tilde{\Lambda}}_{\gamma}
$$
we have
$$
\alpha x\alpha^{-1}=\alpha
z\alpha^{-1}+(\alpha\phi\alpha^{-1})(\alpha
y\alpha^{-1})(\alpha\phi^{-1}\alpha^{-1})\in
{\tilde{\Lambda}}_{\gamma}.
$$

We denote by
$$
{\tilde{\Lambda}}_{\gamma}/(\Phi(k),{\zeta(\gamma)})
$$
the orbit space of ${\tilde{\Lambda}}_{\gamma}/\Phi(k)$ under
conjugation by $\zeta(\gamma)$.

\begin{rem}\label{no-group-remn}
In order to see the difficulties in putting a group structure on
${\tilde{\Lambda}}_{\gamma}/(\Phi(k),{\zeta(\gamma)})$, note that if $(\mu(K),\phi)\in\Phi(k)$ and
$\phi\mu(K)\phi^{-1}\neq\mu(K)\in{\tilde{\Lambda}}_{\gamma}$, then $\mu(K)$ and
$\mu(K)+\phi\mu(K)\phi^{-1}$ are each in the orbit of
$0\in{\tilde{\Lambda}}_{\gamma}/(\Phi(k),{\zeta(\gamma)})$ but $\mu(K)+\mu(K)=2\mu(K)$ may {\em
not} be in the orbit of $0$.
\end{rem}

\subsection{Relative self-linking numbers}\label{rel-self-link-subsec}
We can now define the {\em algebraic relative self-linking numbers} $\mu_{k}(j)$; the adjective
``algebraic'' will usually be omitted in the interest of brevity.
\begin{defi}\label{rel-self-link-defi}
Fix any knot $k$ in $\mathcal{K}_{\gamma}(M)$. For all $j\in\mathcal{K}_{\gamma}(M)$, define
$\mu_{k}(j)$ the {\em relative self-linking number of $j$ with respect to $k$} by
$$
\mu_{k}(j):=\mu(H)\in{\tilde{\Lambda}}_{\gamma}/(\Phi(k),{\zeta(\gamma)})
$$
where $H$ is any singular concordance from $k$ to $j$.
\end{defi}

\subsection{Proof of Theorem~\ref{rel-self-link-thm}}\label{rel-self-link-thm-proof}

\begin{figure}[ht!]
         \centerline{\includegraphics[scale=.45]{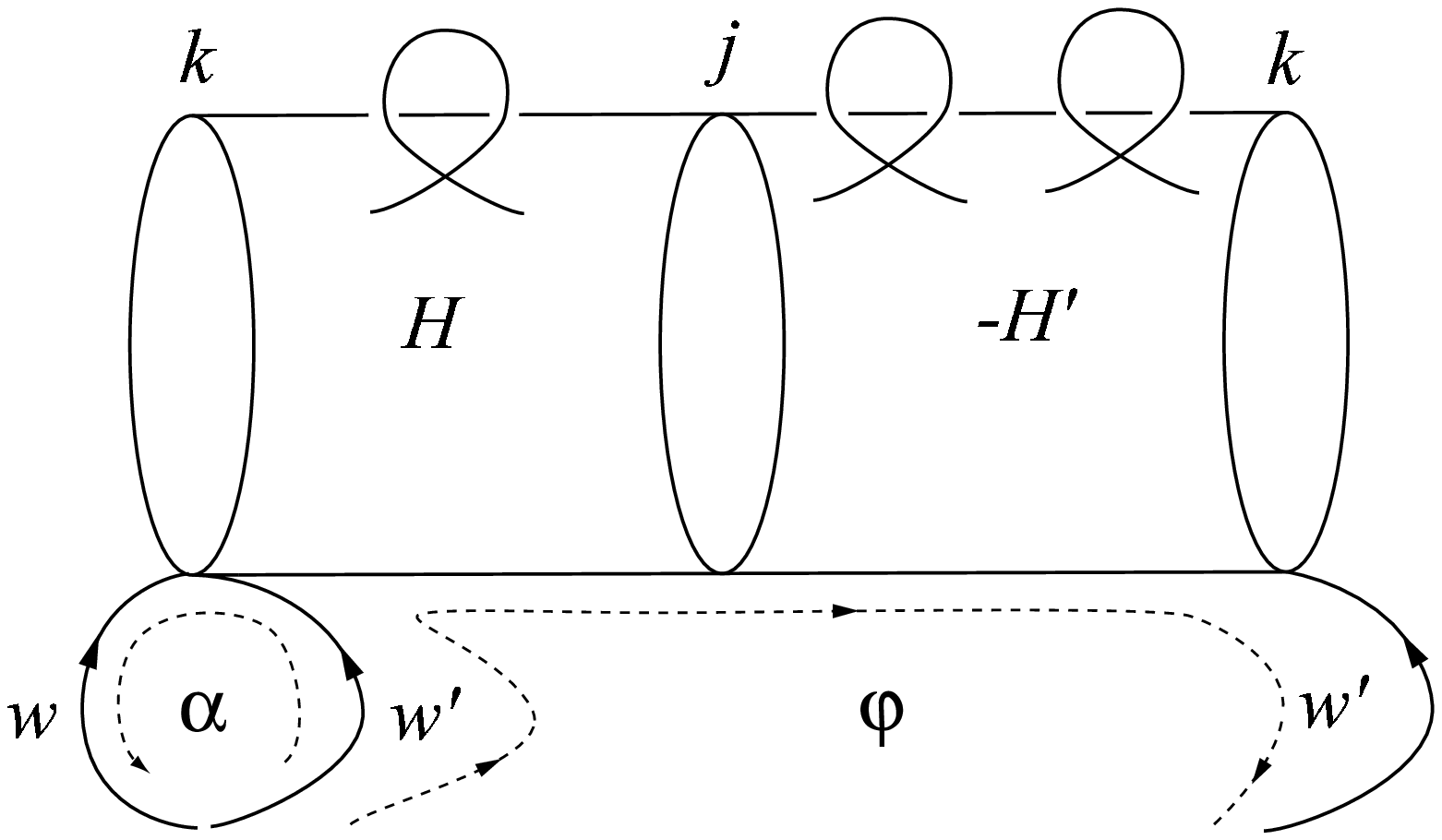}}
         \nocolon \caption{}
         \label{H-conc-fig}

\end{figure}

After unravelling equivalence relations, we have that $\mu_k(j)=\mu_k(j')$ if and only if for any
singular concordances $H$ and $H'$ from $k$ to $j$ and $j'$, respectively, there exist an element
$(z,\phi)\in\Phi(k)$ and an element $\alpha\in\zeta(\gamma)$ such that
$$
\alpha\mu(H)\alpha^{-1}=z+\phi\mu(H')\phi^{-1}\in{\tilde{\Lambda}}_{\gamma}.
$$
Here $\mu(H)$ and $\mu(H')$ are computed using whiskers $w$ and $w'$, respectively, on $k$ such
that $[k]=\gamma$.
\paragraph{Independence of choice of singular concordance}
If $H$ and $H'$ are any two singular concordances from $k$ to $j$, then the composition $K:=H-H'$
is a singular self-concordance of $k$ (see Figure~\ref{H-conc-fig}). Using the whisker $w'$ on $k$
as a whisker for $K$, we have
$$
\mu(K)=\alpha\mu(H)\alpha^{-1}-\phi\mu(H')\phi^{-1}\in {\tilde{\Lambda}}_{\gamma}
$$
where $\alpha\in\zeta(\gamma)$ is determined by the difference of the whiskers $w$ and $w'$ and
$\phi\in\zeta(\gamma)$ is determined by a latitude of $K$ (using the whisker $w'$). Thus, the
relative self-linking number $\mu_k(j)$ does not depend on the choice of singular concordance from
$k$ to $j$ since the equation
$$
\alpha\mu(H)\alpha^{-1}=\mu(K)+\phi\mu(H')\phi^{-1}
$$
in ${\tilde{\Lambda}}_{\gamma}$ shows that
$$
\mu(H)=\mu(H')\in{\tilde{\Lambda}}_{\gamma}/(\Phi(k),{\zeta(\gamma)}).
$$

\paragraph{Concordance invariance}
First we show that concordant knots have the same indeterminacy subgroup. Let $C$ be any
concordance from $k'$ to $k$. Then for any singular self-concordance $K$ of $k$, the composition
$K':=C+K-C$ is a singular self-concordance of $k'$ and, since $\mu(C)=0$, we have
$$
\mu(K')=\mu(C+K-C)=\alpha\mu(K)\alpha^{-1}
$$
where $\alpha\in\zeta(\gamma)$ is determined by a latitude of $C$ which depends on the whiskers for
$k$ and $k'$. By choosing an appropriate whisker for $k'$ it can be arranged that $\alpha=1$, which
means that
$$
(\mu(K'),\phi')=(\mu(K),\phi)\in{\tilde{\Lambda}}_{\gamma}\rtimes \zeta(\gamma)
$$
where $\phi$ and $\phi'$ are determined by latitudes of $K$ and $K'$ respectively. Thus,
$\Phi(k)\leq\Phi(k')$ and by a symmetric argument $\Phi(k)=\Phi(k')$.

Now we also vary $j$ by a concordance. If $C'$ is any concordance from $j$ to $j'$, then for any
singular concordance $H$ from $k$ to $j$ the composition $H':=C+H+C'$ is a singular concordance
from $k'$ to $j'$ and, since $C$ and $C'$ have no singularities, we have
$$
\mu_{k'}(j')=\mu(H')=\alpha\mu(H)\alpha^{-1}=\mu_k(j)
\in{\tilde{\Lambda}}_{\gamma}/(\Phi(k),{\zeta(\gamma)})
$$
where $\alpha\in\zeta(\gamma)$ is determined by a latitude of the concordance $C$ from $k'$ to $k$.
Thus, the induced map in the statement of Theorem~\ref{rel-self-link-thm} is well-defined.

\paragraph{The induced map is onto}
To see that $\mathcal{C}_{\gamma}(M)$ is mapped onto the target space, just apply the same
construction as in the proof Theorem~\ref{absolute-link-thm} (see~\ref{onto-subsubsection}): push
out small arcs of $k$ around loops and create $\pm$-clasps with $k$.\endproof


\section{The indeterminacy sub-group $\Phi(k)$}\label{Phi-sec}
In this section we take up the question of computability of relative self-linking numbers. This
computability depends to large extent on describing the indeterminacy subgroup $\Phi(k)$. The main
goal of the section is Proposition~\ref{Phi-prop} which describes a set of generators for
$\Phi(k)$. The indeterminacies come from two kinds of self-homotopies: {\em Spherical
self-homotopies} (\ref{spherical-self-htpy2}) contribute indeterminacies that can be measured by an
intersection pairing between spheres and knots (\ref{sphere-ints}) and indeterminacies coming from
{\em toriodal self-homotopies} (\ref{essential-self-htpy}) correspond to essential maps of tori and
can be checked on a generating set for the centralizer $\zeta(\gamma)$ of $\gamma=[k]$.
Lemma~\ref{straightening-lemma} and its corollary present the key observation that the spherical
self-homotopies are precisely those with trivial latitudes.

In fact, the generating set of $\Phi(k)$ often simplifies significantly in practice, as described
in the paragraphs following the statement of Proposition~\ref{Phi-prop} and as illustrated in the
examples of Section~\ref{example-sec}.

\subsection{Intersecting knots and spheres}\label{sphere-ints}

Let $k$ be a knot in $\mathcal{K}_{\gamma}(M)$ and $S$ be an immersed 2--sphere in $M$ equipped
with whiskers so that $[k]=\gamma\in\pi_1M$ and $[S]=\sigma\in\pi_2M$. To each intersection point
$p\in k\cap S$, associate an element $g_p\in \pi_1M$ from a loop going from the basepoint of $M$
along the whisker on $S$, along $S$ to $p$ (without changing sheets at any singularities of $S$),
along $k$, then back along the whisker on $k$ to the basepoint of $M$. The sum of the $g_p$, taken
in $\Lambda_{\gamma}$, only depends on the homotopy classes of $k$ and $S$.

\begin{defi}
For $k\in\mathcal{K}_{\gamma}(M)$ and $\sigma \in \pi_2M$, define
$$
\lambda(\sigma,k):=\sum (\sign  p)\cdot g_p \in \Lambda_{\gamma}
$$
where the sum is over all intersection points $p$ between $k$ and any immersed 2--sphere
representing $\sigma$.
\end{defi}
Changing the whisker identifying $[k]=\gamma$ changes
$\lambda(\sigma,k)$ by right multiplication by a difference of
whiskers in $\zeta(\gamma)$.

We will use the following linearity property of
$\lambda(\sigma,k)$:
$$
\lambda(g_1\cdot\sigma_1+g_2\cdot\sigma_2,k)=g_1\lambda(\sigma_1,k)+g_2\lambda(\sigma_2,k)
$$
for all $g_1,g_2\in\pi_1M$ and $\sigma_1,\sigma_2\in\pi_2M$.

\subsection{Spherical self-homotopies}\label{spherical-self-htpy2}
\begin{figure}[ht!]
         \centerline{\includegraphics[scale=.45]{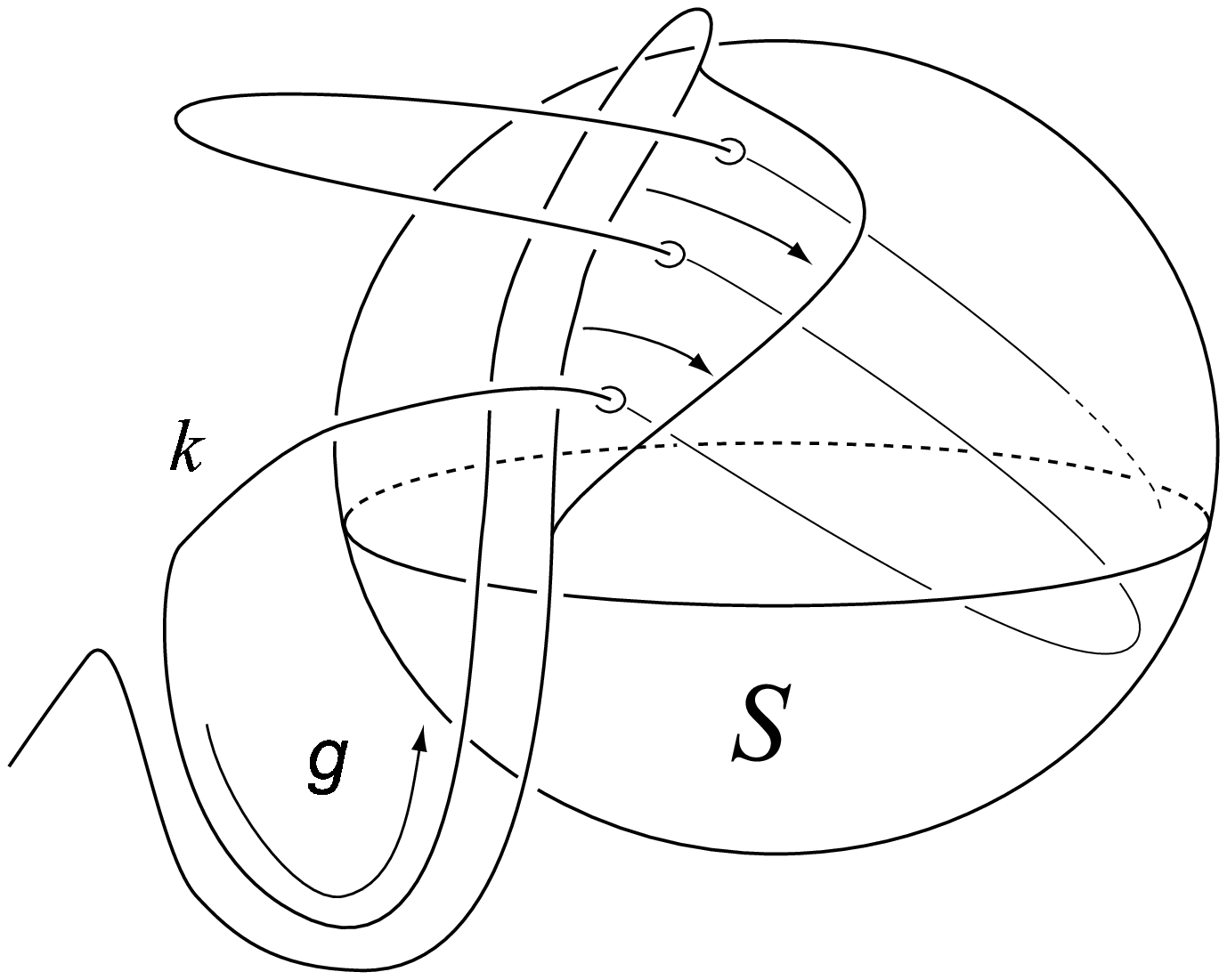}}
         \nocolon \caption{}
         \label{spherical-htpy2-fig}

\end{figure}

Pushing a local strand of a knot $k\in\mathcal{K}_{\gamma}(M)$ around a loop representing
$g\in\pi_1M$, then across an embedded 2--sphere $S$ as indicated in
Figure~\ref{spherical-htpy2-fig}, creates a {\em spherical self-homotopy} of $k$. Denoting by $K$
the trace in $M\times I$ of such a spherical self-homotopy, we have
$$
\mu(K)=\widetilde{\lambda}(g\cdot\sigma,k)\in{\tilde{\Lambda}}_{\gamma}
$$
where $S$ (appropriately based) represents $\sigma\in\pi_2M$ and $\widetilde{\lambda}(\sigma,k)$ is
$\lambda(\sigma,k)$ taken in the quotient ${\tilde{\Lambda}}_{\gamma}$ of $\Lambda_{\gamma}$. Note
that $K$ has a latitude which determines the trivial element of $\pi_1M$. The next lemma and its
corollary give a converse to this observation.

\begin{lem}\label{straightening-lemma}
If a singular self-concordance $K$ of a knot $k\subset M$ has a latitude that represents the
trivial element $1\in\pi_1M$, then $K$ is homotopic (rel $\partial$) to the connected sum of an
embedded 2--sphere and the embedded annulus $A=k \times I$ in $M \times I$.
\end{lem}

\begin{figure}[ht!]
         \centerline{\includegraphics[scale=.35]{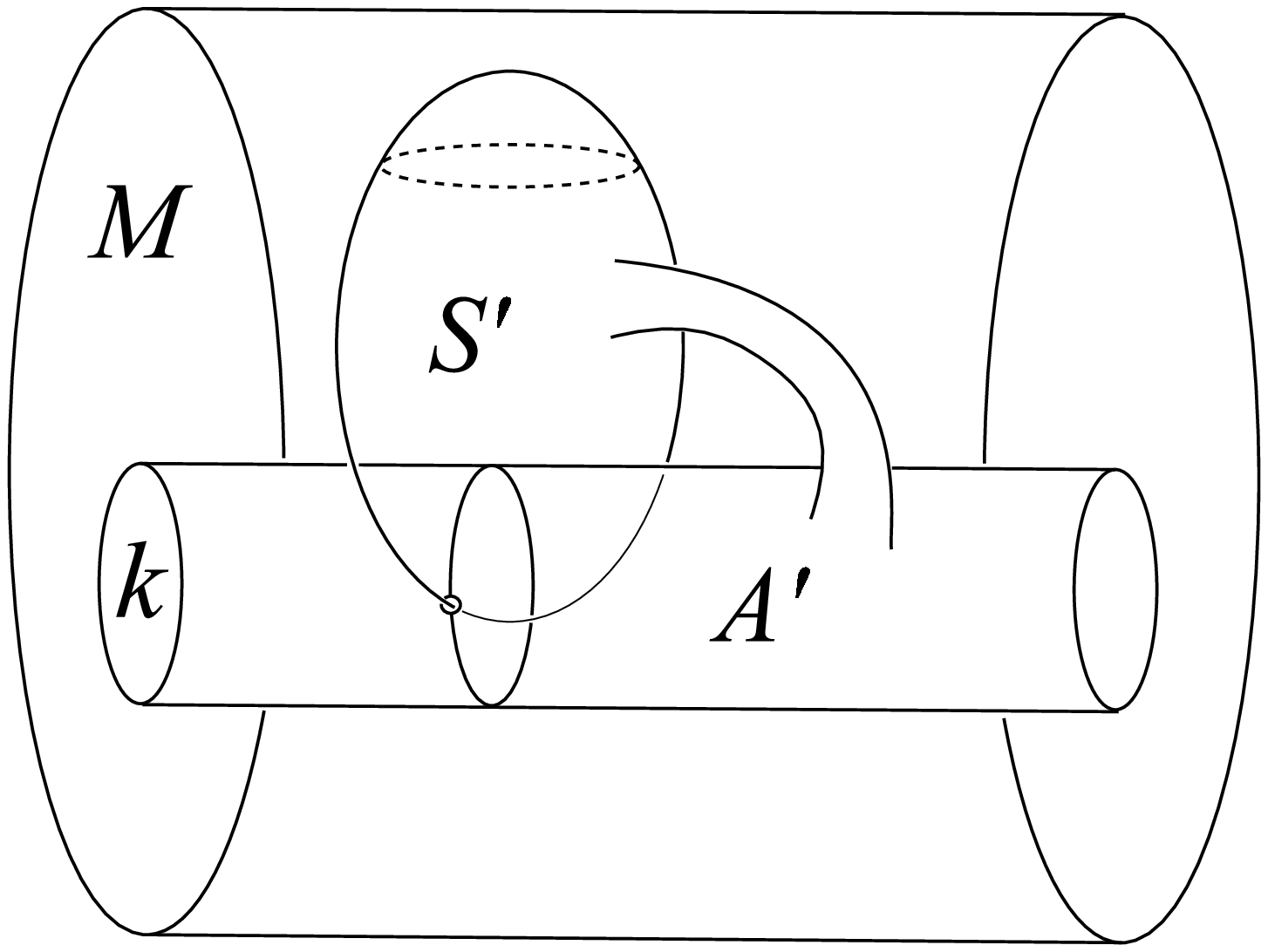}}
         \nocolon \caption{}
         \label{straightening-fig}

\end{figure}

\begin{proof} The projection $T$ of $K$ in $M$ is a map of a torus into
$M$. This map $T$ extends to a map (we will still call $T$) of a torus with a 2--disk $D$ attached
along the homotopically trivial loop which is the projection of the latitude. The obstruction to
further extending $T$ to a map of a solid torus is an immersed 2--sphere $S$ in $M$ which is the
result of using (the image of) $D$ to surger $T$.

If $S$ bounds an immersed 3--ball in $M$, then $T$ extends to a map of a solid torus $S^1 \times
D^2$ with $k$ corresponding to the $S^1$ factor. (This solid torus is the union of the 3--ball and
$D\times I$ joined along $D\times \{0,1\}$.) In this case, by collapsing along the $D^2$ factor (as
indicated in Figure~\ref{disk-sphere-collapse-fig}(a)) we get a homotopy from $T$ to $k$ which
lifts to a homotopy (rel $\partial$) from $K$ to $A$.

\begin{figure}[ht!]
         \centerline{\includegraphics[scale=0.40]{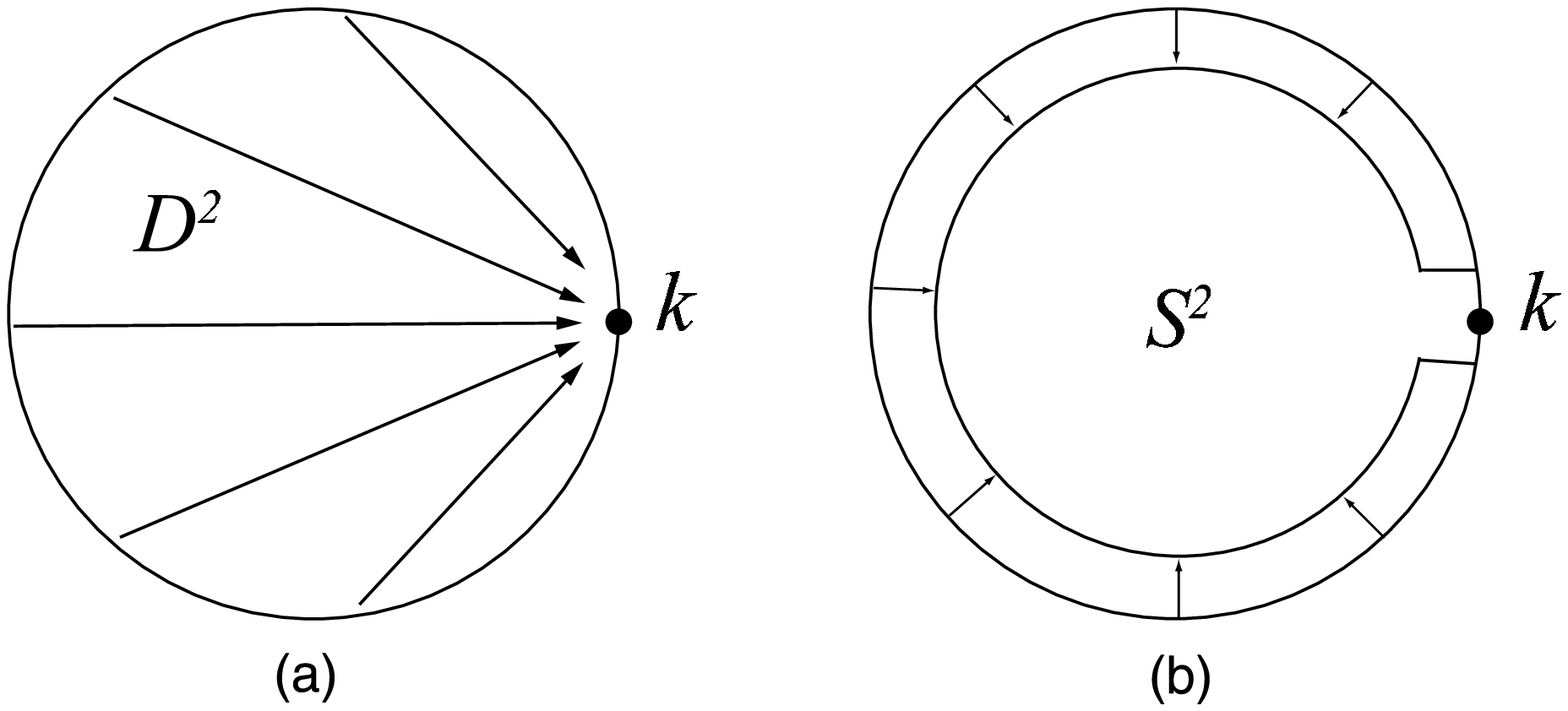}}
         \nocolon \caption{}
         \label{disk-sphere-collapse-fig}

\end{figure}

If $S$ does not bound a map of a 3--ball, then $T$ deformation retracts onto the connected sum of
$S$ with a small band on $k$ (Figure~\ref{disk-sphere-collapse-fig}(b)). In this case the
retraction lifts to a homotopy (rel $\partial$) from $K$ to the connected sum $A'$ of $A$ with a
2--sphere $S'$ that is a lift of $S$ (see Figure~\ref{straightening-fig}). By
Lemma~\ref{embedded-sphere-lemma} we can assume that $S'$ is embedded.
\end{proof}

\begin{cor}\label{straightening-cor}

If $K$ is a singular self-concordance of any knot $k$ in $\mathcal{K}_{\gamma}(M)$ and $K$ has a
latitude which determines the trivial element of $\pi_1M$, then
$\mu(K)=\widetilde{\lambda}(\sigma,k)$ for some $\sigma\in\pi_2M$.

\end{cor}

\begin{proof}
By Lemma~\ref{straightening-lemma} and the homotopy invariance of $\mu$, we have that
$\mu(K)=\mu(A')$, where $A'$ is the connected sum of an embedded 2--sphere $S'$ and the embedded
product annulus $A=k\times I$ in $M\times I$. The embedded 2--sphere $S'$ consists of parallel
copies of embedded generating 2--spheres in different copies of $M$ cross a point in $M\times I$
tubed together as described in the proof of Lemma~\ref{embedded-sphere-lemma}. Thus, the only
singularities of the annulus $A'$ come from intersections between $k$ and the embedded 2--spheres
in $M$. The corollary follows by additivity (\ref{sphere-ints}) with $\sigma=[S']$.
\end{proof}

\subsection{Toroidal self-homotopies}
\label{essential-self-htpy}
 For each non-trivial element $\phi\in\zeta(\gamma)$, the equation
$\gamma\phi\gamma^{-1}\phi^{-1}=1\in\pi_1M$ corresponds (after a choice of a null-homotopy) to a
map of a torus into $M$. For $k\in\mathcal{K}_{\gamma}(M)$, such a map parametrizes a self-homotopy
$K_{\phi}$ which has a latitude representing $\phi$. The next lemma states that the
self-intersection number of the trace of $K_{\phi}$ is unique up to spherical intersections. Such
{\em toroidal self-homotopies} will appear in Proposition~\ref{Phi-prop} below.

\begin{lem}\label{Phi-lem}
For $k\in\mathcal{K}_{\gamma}(M)$, if $K_{\phi}$ and $K'_{\phi}$
are singular self-concordances of $k$ with
$$
\lat[K](w)=\lat[K'](w)\phi\in\zeta(\gamma)
$$
for some whisker $w$ on $k$ then
$$
\mu(K_{\phi})(w)-\mu(K'_{\phi})(w)=\lambda(\sigma,k)
$$
for some $\sigma\in\pi_2M$.
\end{lem}

\begin{proof}
Apply Corollary~\ref{straightening-cor} to $K_{\phi}-K'_{\phi}$.
\end{proof}

\subsection{A generating set for $\Phi(k)$}

\begin{prop}\label{Phi-prop}

Given $k\in\mathcal{K}_{\gamma}(M)$, let $\{\phi_q\}$ be any
generating set for the centralizer $\zeta(\gamma)$ and let $\{K_q
\}$ be any set of traces of self-homotopies of $k$ such that, for
a fixed whisker $w_0$ on $k$, $\lat[K_q](w_0)=\phi_q$ for each
$q$. Let $\{\sigma_s\}$ be any generating set for $\pi_2M$ as a
module over $\pi_1M$. Denote by $\mathcal{S}_0$ the set
$$
\mathcal{S}_0:=\{(\mu(K_q)(w_0),\phi_q),(g\widetilde{\lambda}(\sigma_s,k)(w_0),1)
\}
$$
where $g$ ranges over $\pi_1M$.

Then $\Phi(k)$ is generated by
$$
\bigcup_{\alpha\in\zeta(\gamma)}\alpha\mathcal{S}_0\alpha^{-1}
$$
where the conjugation $\alpha\mathcal{S}_0\alpha^{-1}$ of
$\mathcal{S}_0$ by $\alpha$ denotes conjugation of each element of
$\mathcal{S}_0$ by $\alpha$.

\end{prop}

The generating set described in Proposition~\ref{Phi-prop} will
often simplify significantly in practice. Usually $\mathcal{S}_0$
will contain either non-trivial spherical indeterminacies or
non-trivial toroidal indeterminacies but not both. Also, the
conjugation action by $\zeta(\gamma)$ is often trivial; for
instance if $\gamma$ is primitive and $\zeta(\gamma)$ is cyclic.

If $M$ is irreducible, then $\zeta(\gamma)$ is finitely generated \cite{Ja} and there are no
spherical indeterminacies so $\mathcal{S}_0$ is a finite set. However, if $M$ is not irreducible
and $k$ has non-trivial intersections with a 2--sphere, then $\Phi(k)$ is, in general, infinitely
presented as will be illustrated in Section~\ref{example-sec}.

 Before giving a proof of
Proposition~\ref{Phi-prop}, we give an example which will also serve to illustrate constructions
used later in the proof of Theorem~\ref{M-thm}.

\begin{example}\label{rel-example1-trivial}
Consider the knot $k\in\mathcal{K}_{xyz}(M)$ in Figure~\ref{rel-example1-fig} of the introduction,
where $M$ is the product of a thrice punctured 2--disk with the circle. In this case $M$ is
irreducible, so there are no generators of $\Phi(k)$ coming from spherical intersections and
$\zeta(xyz)$ is isomorphic to $\Z\oplus\Z$ generated by $xyz$ and $t$; $\Phi(k)$ can be computed by
constructing just two self-homotopies of $k$ with latitudes representing $xyz$ and $t$
respectively. Two such self-{\em isotopies} of $k$ can easily be constructed, one by isotoping $k$
along itself longitudinally and the other by isotoping $k$ around the circle fiber of $M$. Since
the trace of a self-isotopy is an embedded annulus, and hence has vanishing self-intersection
number in $\widetilde{\Lambda}_{\gamma}$, the set $\mathcal{S}_0$ of Proposition~\ref{Phi-prop}
contains only the two elements $(0,xzy)$ and $(0,t)$ which shows that the action $\Phi(k)$ reduces
to conjugation by $\zeta(xyz)$. Since $t$ is central in $\pi_1M$ and $xyz$ is primitive, it follows
that the conjugation action of $\zeta(xyz)$ on $\widetilde{\Lambda}_{xyz}$ is trivial, justifying
the statement that the target space for relative self-linking numbers with respect to $k$ is
$\widetilde{\Lambda}_{xyz}$ itself.
\end{example}

\subsection{Proof of Proposition~\ref{Phi-prop}}

Let $(\mu(K)(w),\phi)$ be any generator of
$\Phi(k)\leq{\tilde{\Lambda}}_{\gamma}\rtimes \zeta(\gamma)$; that
is, $K$ is a singular self-concordance of $k$ such that
$\lat[K](w)=\phi\in\zeta(\gamma)$ for some whisker $w$.

First consider the case when $w=w_0$: Since the $\phi_q$ generate $\zeta(\gamma)$, there exists a
finite sequence $K_{q_i}$ of the $\{K_q \}$ such that
$$
K':=K+\sum_iK_{q_i}
$$
is a singular self-concordance of $k$ with
$\lat[K'](w_0)=1\in\pi_1M$.

By Corollary~\ref{straightening-cor},
$$
\mu(K')(w_0)=\widetilde{\lambda}(\sigma,k)(w_0)=\sum_s\widetilde{\lambda}(g_s\cdot
\sigma_s,k)(w_0)
$$
for some $\sigma\in\pi_2M$.

We have the following equation in
${\tilde{\Lambda}}_{\gamma}\rtimes \zeta(\gamma)$:
$$
(\mu(K)(w_0),\phi)\prod_i(\mu(K_{q_i})(w_0),\phi_{q_i})=\prod_s(g_s\widetilde{\lambda}(\sigma_s,k)(w_0),1)
$$
which can be solved for $(\mu(K)(w_0),\phi)$ in terms of elements
(and their inverses) from $\mathcal{S}_0$.

For the general case where $w$ and $w_0$ differ by $\alpha\in\zeta(\gamma)$, the analogous
construction expresses $(\mu(K)(w),\phi)$ in terms of elements from
$\alpha\mathcal{S}_0\alpha^{-1}$. \endproof

\section{Spherical knots and the proof of Theorem~\ref{spherical-thm}}\label{spherical-sec}
The proof of Theorem~\ref{spherical-thm} is given in this section after first defining {\em
spherical knots} and stating a lemma which isolates the properties of spherical knots to be used in
the proof. The proof of the lemma is deferred until the end of the section. Here is the key point:
For a spherical knot, every element of the centralizer can be realized as a latitude of a
self-homotopy with {\em vanishing} self-intersection number; this means that the conjugation
effect, which in general obstructs additivity of $\mu$ under composition of singular concordances
(\ref{composing-concs}), can always be ``cancelled out'' by inserting these self-homotopies
accordingly.
\subsection{Spherical knots}\label{spherical-subsec}
\begin{defi}\label{spherical-defi}
A knot $k_0\in \mathcal{K}_{\gamma}$ is {\em spherical} if there exists a whisker $w_0$ on $k_0$
such that for any $\phi\in\zeta(\gamma)$ there exists a singular self-concordance $K^0_{\phi}$ of
$k_0$ with $\lat[K^0_{\phi}](w_0)=\phi$ and $\mu(K^0_{\phi})(w_0)=0\in{\tilde{\Lambda}}_{\gamma}$.
\end{defi}
By Proposition~\ref{Phi-prop}, ``$k_0$ is spherical'' means that ``$\Phi(k_0)$ contains only
spherical indeterminacies,'' that is, indeterminacies due to intersections between $k_0$ and
2--spheres in $M$.
\begin{example}\label{rel-example1-spherical}
The knot $k$ in Figure~\ref{rel-example1-fig} of the introduction is spherical by the discussion in
Example~\ref{rel-example1-trivial} and any knot in $S^1\times S^2$ is spherical as will be
discussed in \ref{S^1-S^2}.
\end{example}

\paragraph{Remark} For any two spherical knots $k_0$ and $k_0'$ in $\mathcal{K}_{\gamma}$, we have
$\Phi(k_0)=\Phi(k_0')$ since intersections between knots and spheres only depend on the homotopy
classes.

The next lemma shows that, for spherical knots, all elements of the indeterminacy sub-group can be
realized by actual singular self-concordances. (This is {\em not} true in general for non-spherical
knots! See remarks \ref{no-group-remn} and \ref{spherical-realization-remn}.)

\begin{lem}\label{spherical-lem}
If $k\in\mathcal{K}_{\gamma}(M)$ is spherical, then
\begin{enumerate}
  \item for any $\alpha\in\zeta(\gamma)$ and any whisker $w$ on
  $k$, there exists a singular self-concordance $K^0_{\alpha}$ of
  $k$ such that $\lat[K^0_{\alpha}](w)=\alpha$ and
  $$
  \mu(K^0_{\alpha})(w)=0\in\tilde{\Lambda}_{\gamma}.
  $$
  \item for every $(z,\phi)\in\Phi(k)$ and every whisker $w$ on
  $k$, there exists a singular self-concordance $K_{\phi}^z$ of $k$
  such that $\lat[K_{\phi}^z](w)=\phi$ and
  $$
  \mu(K_{\phi}^z)(w)=z\in\tilde{\Lambda}_{\gamma}.
  $$

\end{enumerate}
\end{lem}
The proof of this lemma is deferred until the end of this section.

\subsection{Proof of Theorem~\ref{spherical-thm}}

\paragraph{Proof of (i) $\mathbf{\Rightarrow}$ (ii)} We will construct a singular concordance
$H^0$ between $j$ and $j'$ such that $\mu(H^0)=0\in\tilde{\Lambda}_{\gamma}$ which by
Proposition~\ref{annuli-selfint-prop} implies (ii).

The equality of $\mu_k(j)$ and $\mu_k(j')$ means that for any singular concordances $H$ and $H'$
from $k$ to $j$ and $j'$, respectively, there exists $\alpha\in\zeta(\gamma)$ and
$(z,\phi)\in\Phi(k)$ such that
$$
\alpha\mu(H)\alpha^{-1}=(z,\phi)\cdot\mu(H')\in{\tilde{\Lambda}}_{\gamma}.
$$
Since this holds for {\em any} such $H$ and $H'$, we may assume that $H$ and $H'$ have the same
whisker $w$ on $k$ so that
$$
\alpha(\mu(H)(w))\alpha^{-1}=z+\phi(\mu(H')(w)\phi^{-1}\in{\tilde{\Lambda}}_{\gamma}.
$$

By Lemma~\ref{spherical-lem}, there exist singular self-concordances $K_{\alpha}^0$ and
$K_{\phi}^z$ of $k$  such that $\mu(K_{\alpha}^0)(w)=0$ and $\mu(K_{\phi}^z)(w)=z$ in
${\tilde{\Lambda}}_{\gamma}$ with $\lat[K_{\alpha}^0](w)=\alpha$ and $\lat[K_{\phi}^z](w)=\phi$.

\begin{figure}[ht!]
         \centerline{\includegraphics[scale=0.50]{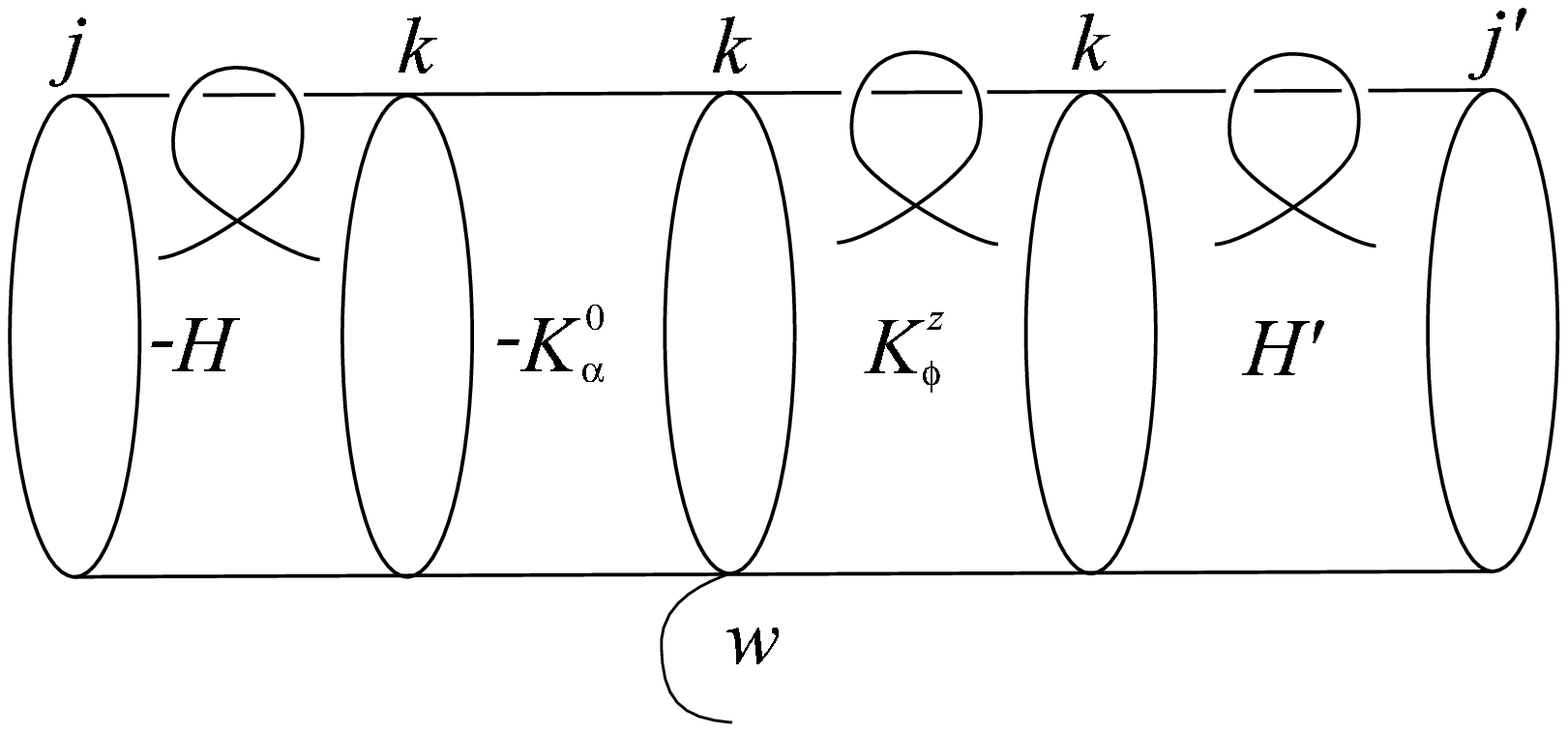}}
         \nocolon \caption{}
         \label{sing-conc-comp-fig}

\end{figure}

The composition $H^0:=-H-K_{\alpha}^0+K_{\phi}^z+H'$ is a singular concordance from $j$ to $j'$
(see Figure~\ref{sing-conc-comp-fig}). Using the whisker $w$ in the intersection
$(-K_{\alpha}^0\cap K_{\phi}^z)\subset H^0$ as a whisker for $H^0$, we have
$$
\begin{array}{rcl}
  \mu(H^0)(w) & = & \mu(-H-K_{\alpha}^0+K_{\phi}^z+H')(w) \\
   & = & -\alpha(\mu(H)(w))\alpha^{-1}+0+z+\phi(\mu(H')(w))\phi^{-1} \\
   & = & 0\in{\tilde{\Lambda}}_{\gamma}
\end{array}
$$
which, by Proposition~\ref{annuli-selfint-prop}, completes the proof that (i) implies (ii).\endproof

\paragraph{Proof of (ii) $\mathbf{\Rightarrow}$ (i)} Let $J$ be a singular concordance from $j$ to
$j'$ with all singularities of $J$ paired by Whitney disks. By
Proposition~\ref{annuli-selfint-prop}, $\mu(J)=0\in{\tilde{\Lambda}}_{\gamma}$. If $H$ is any
singular concordance from $k$ to $j$, the composition $H+J$ is a singular concordance from $k$ to
$j'$ and we have
$$
\mu_k(j)=\mu(H)=\mu(H+J)=\mu_k(j')
$$
since all singularities of $J$ come in cancelling
pairs.\endproof

\subsection{Proof of Lemma~\ref{spherical-lem}}

\paragraph{Proof of property (i) of Lemma~\ref{spherical-lem}} Property (i) follows from the
fact that changing whiskers on $k$ corresponds to conjugation by $\zeta(\gamma)$ which induces
isomorphisms of ${\tilde{\Lambda}}_{\gamma}$ and $\zeta(\gamma)$.

Specifically, given $\alpha\in\zeta(\gamma)$ and any whisker $w$ on $k$, if $\psi\in\zeta(\gamma)$
is determined by the difference $[w-w_0]$ of $w$ and $w_0$ (where $w_0$ is the whisker of
Definition~\ref{spherical-defi}), then, for $K^0_{\psi^{-1}\alpha\psi}$ with
$\lat[K^0_{\psi^{-1}\alpha\psi}](w_0)=\psi^{-1}\alpha\psi$ and
$\mu(K^0_{\psi^{-1}\alpha\psi})(w_0)=0\in{\tilde{\Lambda}}_{\gamma}$ as guaranteed by
Definition~\ref{spherical-defi}, we have
$$
\lat[K^0_{\psi^{-1}\alpha\psi}](w)=\psi\lat[K^0_{\psi^{-1}\alpha\psi}](w_0)\psi^{-1}=\alpha
$$
and
$$
\mu(K^0_{\psi^{-1}\alpha\psi})(w)=\psi\mu(K^0_{\psi^{-1}\alpha\psi})(w_0)\psi^{-1}=0
\in{\tilde{\Lambda}}_{\gamma}.\eqno{\qed}
$$

\paragraph{Proof of property (ii) of Lemma~\ref{spherical-lem}}\label{spherical-(ii)}

Let $K^z_{\phi}$ be a singular self-concord\-ance of $k$ with
$$
(z,\phi)=(\mu(K^z_{\phi})(w_0),\lat[K^z_{\phi}](w_0))
\eqno(\mbox{*})
$$
for some whisker $w_0$ on $k$ and let $w$ be any other whisker on $k$ with
$\psi=[w-w_0]\in\zeta(\gamma)$. By property (i) of Lemma~\ref{spherical-lem}, there exists a
singular self-concordance $K^0_{\psi^{-1}}$ of $k$ such that $\lat[K^0_{\psi^{-1}}](w)=\psi^{-1}$
and $\mu(K^0_{\psi^{-1}})(w)=0\in{\tilde{\Lambda}}_{\gamma}$. The composition
$K':=K^0_{\psi^{-1}}+K^z_{\phi}-K^0_{\psi^{-1}}$ satisfies
$$
\lat[K'](w)=\psi^{-1}(\lat[K^z_{\phi}](w))\psi=\psi^{-1}\psi(\lat[K^z_{\phi}](w_0))\psi^{-1}\psi=\phi
$$
and
$$
\mu(K')(w)=\psi^{-1}(\mu(K^z_{\phi})(w))\psi=\psi^{-1}\psi(\mu(K^z_{\phi})(w_0))\psi^{-1}\psi=z.
$$
Thus, to prove property (ii) of Lemma~\ref{spherical-lem} it suffices to show that it holds for a
fixed whisker $w_0$, that is, it suffices to show that for any $(z,\phi)\in\Phi(k)$ there exists
$K^z_{\phi}$ satisfying the above equation (*) for $w_0$.

Since $k$ is spherical, by Proposition~\ref{Phi-prop} $\Phi(k)$ is
generated by the set
$$
\mathcal{S}= \{(0,\phi_q), (\alpha
g\widetilde{\lambda}(\sigma_s,k)(w_0)\alpha^{-1},1)\}
$$
for some whisker $w_0$ on $k$, where $\alpha$ ranges over $\zeta(\gamma)$, the $\phi_q$ generate
$\zeta(\gamma)$, the $\sigma_s$ generate $\pi_2M$ as a module over $\pi_1M$, and $g$ ranges over
$\pi_1M$.

The cases $(z,\phi)=(0,\phi_q)$ (or more generally $z=(0,\phi)$ for any $\phi\in\zeta(\gamma)$ and
any $w$) are covered by property (i) of Lemma~\ref{spherical-lem} (which is essentially built into
Definition~\ref{spherical-defi}).

If $(z,\phi)=(g\widetilde{\lambda}(\sigma_s,k)(w_0),1)$ then, as described
in~\ref{spherical-self-htpy2}, there exists a spherical self-homotopy of $k$ whose trace $K_{gS_s}$
satisfies $\lat[K_{gS_s}](w_0)=1$ and
$\mu(K_{gS_s})(w_0)=g\widetilde{\lambda}(\sigma_s,k)(w_0)\in\tilde{\Lambda}_{\gamma}$ where $gS_s$
is a 2--sphere in $M$ representing $g\cdot\sigma_s\in\pi_2M$.

Conjugating by a singular self-concordance $K^0_{\alpha}$ of $k$ from property (i) yields a
singular self-concordance $K^0_{\alpha}+K_{gS_i}-K^0_{\alpha}$ which satisfies (*) for
$(z,\phi)=(\alpha g \widetilde{\lambda}(\sigma_s,k)(w_0)\alpha^{-1},1)$.

We have found singular self-concordances satisfying (*) for all elements of the generating set
$\mathcal{S}$ for $\Phi(k)$. It only remains to observe that if $K$ and $K'$ are any singular
self-concordances of $k$ with $\lat[K](w_0)=\phi$ and $\lat[K'](w_0)=\phi'$, then $K+K'$ satisfies
(*) for
$$
(z,\phi)=(\mu(K)(w_0),\phi)(\mu(K')(w_0),\phi')=(\mu(K)(w_0)+\phi(\mu(K')(w_0))\phi^{-1},\phi\phi').\eqno{\qed}
$$

\begin{rem}\label{spherical-realization-remn}
It follows from Lemma~\ref{spherical-lem}, together with the observation that
$(\mu(-K)(w),\lat[-K](w))$ is an inverse for $(\mu(K)(w),\lat[K](w))$, that if $k$ is spherical,
then the subset
$$
\{(\mu(K)(w),\lat[K](w))\}\subset{\tilde{\Lambda}}_{\gamma}\rtimes \zeta(\gamma),
$$
as $K$ ranges over all singular self-concordances of $k$ in $M$ and $w$ ranges over all whiskers
identifying $[k]=\gamma\in\pi_1M$, is actually {\em equal} to the subgroup $\Phi(k)$. The subtlety
here is that if a knot is {\em not} spherical, then its indeterminacy subgroup can contain elements
that are not realized by actual singular self-concordances of the knot.
\end{rem}


\section{The class $\mathcal{M}$ and proof of Theorem~\ref{M-thm}}\label{M-sec}
The proof of Theorem~\ref{M-thm} is given in this section after the definition of the set
$\mathcal{M}$. The subset of {\em irreducible} manifolds in $\mathcal{M}$ was denoted $\mathcal{N}$
in \cite{KL1} and \cite{KL2}.

\subsection{The class $\mathcal{M}$}

\begin{defi}\label{M-defi}

Let $\mathcal{M}$ denote the set of all 3--manifolds $M$ such that

\begin{enumerate}

\item  $M$ does not contain any circle bundle over a non-orientable
surface whose total space is orientable,

\item  $M$ does not contain any of the Seifert manifolds
$M(S^2,(3,1),(3,1),$\break $(3,-2))$, $M(S^2,(2,1),(4,-1),(4,-1))$ or
$M(S^2,(2,1),(3,-1),(6,-1))$.

\end{enumerate}

\end{defi}
Here ``contain'' means either as a submanifold or as a (punctured) prime factor and the notation
$M(S^2,(3,1),(3,1),(3,-2))$, for instance, refers to the Seifert fibered manifold with orbit surface $S^2$ and
three singular fibers having Seifert invariants $(3,1)$, $(3,1)$ and $(3,-2)$. The excluded manifolds are exactly
those which contain essential tori that contribute non-trivial indeterminacies to $\Phi(k)$ for all
$k\in\mathcal{K}_{\gamma}(M)$ for certain $\gamma$ (see \cite{KL1, KL2, S2}).

\subsection{Proof of Theorem~\ref{M-thm}}

To show that a knot $k_0$ is spherical it is sufficient, by
Proposition~\ref{Phi-prop}, to find self-homotopies of $k_0$ whose
traces $K_q^0$ satisfy
$$
(\mu(K_q^0)(w_0),\lat[K_q^0](w_0))=(0,\phi_q)\in\tilde{\Lambda}_{\gamma}\rtimes
\zeta(\gamma)
$$
for some whisker $w_0$ on $k_0$ where the $\phi_q$ generate $\zeta(\gamma)$. We will in fact be
able to find self-{\em isotopies} of $k_0$ in all cases.

The proof breaks into different cases according to the structure
of $\zeta(\gamma)$.

\paragraph{Null-homotopic knots} For the case $\gamma=1$ take
$k_0$ to be the unknot, that is, $k_0$ bounds an embedded 2--disk in $M$ (all such $k_0$ are
clearly isotopic in any connected manifold). Now $\zeta(\gamma)=\pi_1M$ and $k_0$ can be isotoped
into a neighborhood of a point and then isotoped around any loop which represents a generator of
$\pi_1M$, completing the null-homotopic case.

\paragraph{Remark} If $k_0$ is the unknot, then the action of $\Phi(k_0)$ reduces to conjugation (since
null-homotopic knots have trivial intersections with spheres) so relative self-linking numbers
reduce to (absolute) self-linking numbers: $\mu_{k_0}(j)=\mu(j)$ for all $j\in\mathcal{K}_{1}(M)$.

From now on in the proof we will assume that $\gamma\neq 1$.

\paragraph{Cyclic centralizers and longitudinal
self-isotopies}\label{longitudinal-self-isotopies}

Consider now the case where $\zeta(\gamma)$ is a cyclic subgroup of $\pi_1(M)$ generated by an
element $\rho$ (possibly $\rho=\gamma$). Take $k_0\in\mathcal{K}_{\gamma}(M)$ to be contained in a
tubular neighborhood $R$ of any knot representing $\rho$. An ambient self-isotopy of $R$ in $M$
that fixes $R$ set-wise but moves once around in the direction of $\rho$ induces a self-isotopy of
$k_0$ whose trace $K^0_{\rho}$ has a latitude representing $\rho$. Such a self-isotopy will be
referred to as a {\em longitudinal self-isotopy}. This completes the case where $\zeta(\gamma)$ is
cyclic.

\begin{rem}\label{longitudinal-rem}

Note that in the case where $\zeta(\gamma)$ is cyclic there exists
a spherical knot in $\mathcal{K}_{\gamma}(M)$ for {\em any} $M$.

\end{rem}

So we may now assume that $\zeta(\gamma)$ and $\pi_1(M)$ are not
cyclic groups.

\paragraph{The non-cyclic centralizer case}

In this setting, the classical 3--manifold structure theorems imply that $\zeta(\gamma)$ is carried
by an embedded Seifert fibered submanifold of $M$: By the Kneser/Milnor prime factorization
theorem, $M$ has a unique connected sum decomposition into irreducible submanifolds and copies of
$S^1\times S^2$ which give a decomposition of $\pi_1M$ as a free product. If $\gamma$ projected
non-trivially to more than one free factor, then $\zeta(\gamma)$ would be cyclic (by
Corollary~4.1.6 of \cite{MKS}), so we may assume that $\gamma$ is carried by an irreducible factor
$M_{\gamma}\subset M$. The irreducible submanifold $M_{\gamma}$ carries $\zeta(\gamma)$, and if
$M_{\gamma}$ is Haken, then by the work of Jaco-Shalen and Johansson $\zeta(\gamma)$ is carried by
a Seifert fibered submanifold of $M_{\gamma}$ (e.g.\ Theorem 4.1 of \cite{JS2}). If $M_{\gamma}$ is
not Haken, then by the theorems of Casson-Jungreis \cite{CJ} and Gabai \cite{G} either
$\zeta(\gamma)$ is carried by a collar of a torus boundary component of $M_{\gamma}$ (which is
Seifert fibered as $S^1\times I \times S^1$) or $M_{\gamma}$ is itself Seifert fibered.

Let $N\subset M_{\gamma}\subset M$ denote the embedded Seifert
fibered sub-manifold that carries $\zeta(\gamma)$. The orbit
surface $F$ of $N$ must be orientable since $M$ is in
$\mathcal{M}$. Let $N_0$ denote the complement of tubular
neighborhoods of the singular fibers (or one regular fiber if
there are no singular ones) of $N$ and let $F_0$ denote $F\cap
N_0$. Then $N_0$ is the trivial circle bundle over $F_0$ and
generators for $\pi_1N$ are represented by a circle fiber and
loops representing generators of $\pi_1F_0$. Denote by
$t\in\pi_1N$ the element represented by a regular fiber (a circle
fiber of $N_0$).

We may assume that $\gamma$ is not represented by a singular fiber
because such elements of $\pi_1(M)$ have cyclic centralizers.

Consider the case where $\gamma$ is in the cyclic subgroup $\langle t \rangle$ generated by a
regular fiber. Since $F$ is orientable, $\langle t \rangle$ is central in $\pi_1N$, so in this case
$\zeta(\gamma)=\pi_1N$. Take $k_0\in\mathcal{K}_{\gamma}(M)$ to be contained in a neighborhood of a
regular fiber. As in the cyclic centralizer case discussed above, the trace $K_t^0$ of a
longitudinal self-isotopy of $k_0$ around the fiber direction of $N_0$ has a latitude representing
$t$. Any generator $\phi$ of $\pi_1N$ that is carried by $F_0$ can be represented by a latitude of
the trace $K_{\phi}^0$ of a self-isotopy of $k_0$ that pushes the neighborhood of the fiber
containing $k_0$ around $N_0$ over any loop in $F_0$ representing $\phi$.

Now consider the case where $\gamma$ is not carried by any fiber. Then $\zeta(\gamma)$ is
isomorphic to $\Z\oplus\Z$ unless $F$ is a closed torus so that $N=T^3$ (possibly punctured); this
is because $M$ contains no non-vertical tori (and by our assumption that $\pi_1M$ is not cyclic).
The theorem holds for $M=T^3$, since any knot in $T^3$ can be self-isotoped around the circle
factors which generate $\pi_1(T^3)\cong \Z\oplus\Z\oplus\Z$, so we will assume that $F$ is not a
closed torus. One of the factors of $\zeta(\gamma)$ is generated by $t$ and the other is generated
by a root $\rho$ of $\gamma$. Take $k_0\in\mathcal{K}_{\gamma}(M)$ to be contained in a tubular
neighborhood of any knot in $N_0$ that represents $\rho$. Then the trace $K_t^0$ of a self-isotopy
of $k_0$ around the fiber has a latitude representing $t$ and the trace $K_{\rho}^0$ of a
longitudinal self-isotopy of $k_0$ has a latitude representing $\rho$. \endproof


\section{Examples}\label{example-sec}

\subsection{Knots in $\mathbf{S^1\times S^2}$}\label{S^1-S^2}

\begin{figure}[ht!]
         \centerline{\includegraphics[scale=.50]{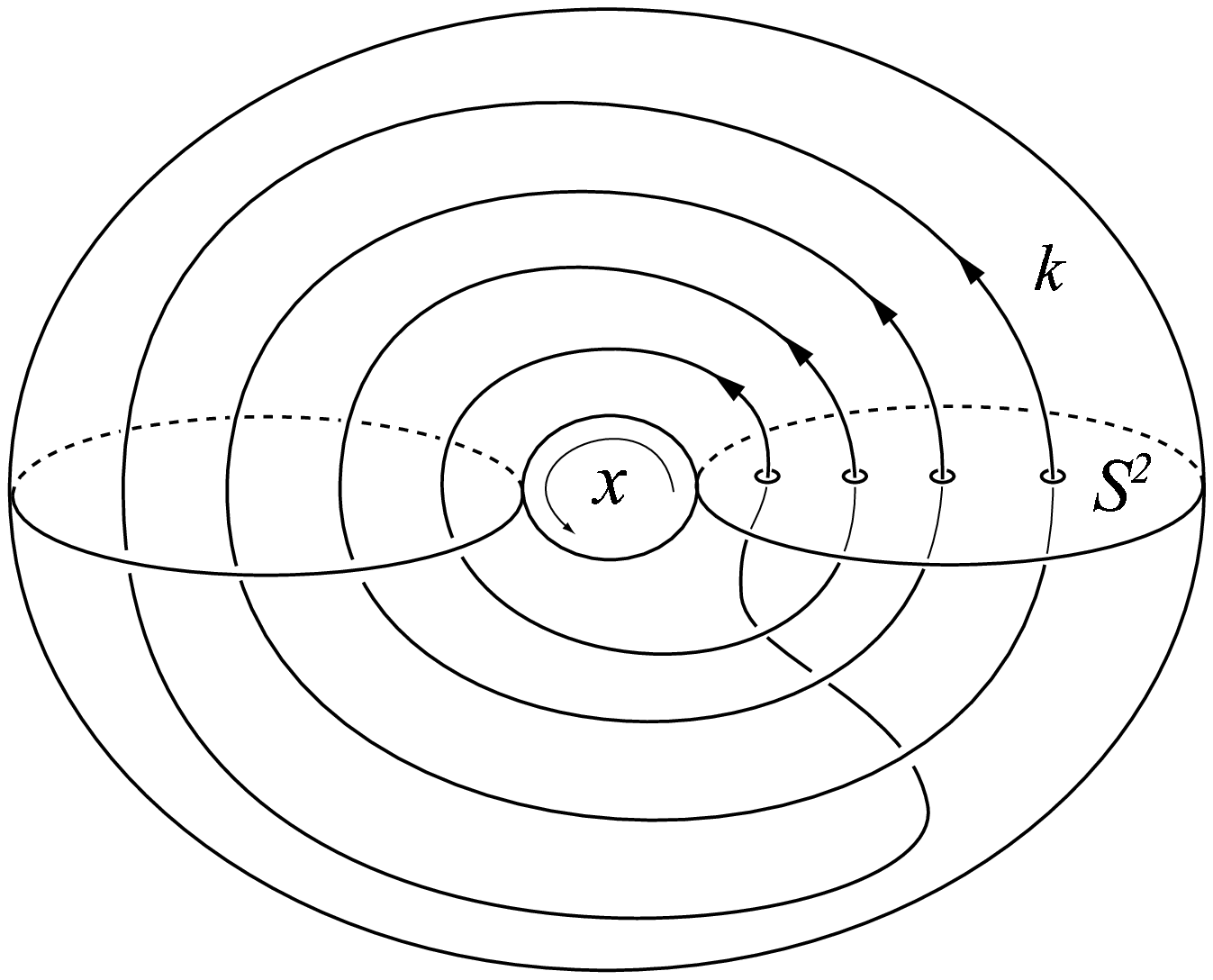}}
         \nocolon \caption{}
         \label{S1crossS2-fig}

\end{figure}

Consider $M=S^1\times S^2\in\mathcal{M}$ with $\pi_1(M)$ generated (multiplicatively) by $x=[S^1]$
and $\pi_2(M)$ generated by $\sigma=[S^2]$. If $\gamma=x^n$, then $\langle x^n \rangle$ is normal
in $\pi_1(M)$ and ${\tilde{\Lambda}}_{x^n}$ is generated by coset representatives
$\{x,x^2,\ldots,x^{n/2}\}$ if $n$ is even, or $\{x,x^2,\ldots,x^{(n-1)/2}\}$ if $n$ is odd (recall
the relations $1=0$ and $g=g^{-1}$). Since all knots in $S^1\times S^2$ can be isotoped around the
$S^1$ direction using the product structure, any $k$ in $\mathcal{K}_{x^n}(S^1\times S^2)$ is
spherical. Since $\pi_1M$ is abelian, conjugation actions are trivial so
$({\tilde{\Lambda}}_{x^n}/(\Phi(k),\zeta(x^n))={\tilde{\Lambda}}_{x^n}/\Phi(k)$. By
Proposition~\ref{Phi-prop}, we see that the generators of $\Phi(k)$ are
$(x^m\lambda(\sigma,k),1)=(x^m(1+x+x^2+\cdots+x^{n-1}),1)$ as can be computed from
Figure~\ref{S1crossS2-fig} (showing a knot $k$ in $\mathcal{K}_{x^4}(S^1\times S^2)$). So in this
case the target orbit space ${\tilde{\Lambda}}_{x^n}/\Phi(k)$ is a group, the quotient of the free
abelian group $\Z[x,x^2,\ldots,x^{n/2}]$ or $\Z[x,x^2,\ldots,x^{(n-1)/2}]$ by a single relation
$x^m(1+x+x^2+\cdots+x^{n-1})= x+x^2+\cdots+x^{n-1}=0$. (This computation can be used to affirm the
conjecture at the end of \cite{KL1}.) Note that if $n$ is odd, then this relation introduces a
2-torsion element:
$$
x+x^2+\cdots+x^{n-1}=2(x+x^2\cdots+x^{(n-1)/2})=0
$$
and if $n$ is even, the relation eliminates the middle degree generator:
$$
x^{n/2}=-2(x+x^2+\cdots+x^{(n/2)-1}).
$$

For instance, if $n=1$ or $n=2$, then the target space is trivial (and any two knots in
$\mathcal{K}_{x}(S^1\times S^2)$ or in $\mathcal{K}_{x^2}(S^1\times S^2)$ are $W\!$--equivalent).

If $n=3$, then the target space is $\Z_2$ since $2x=0$.

\begin{figure}[ht!]
         \centerline{\includegraphics[scale=.60]{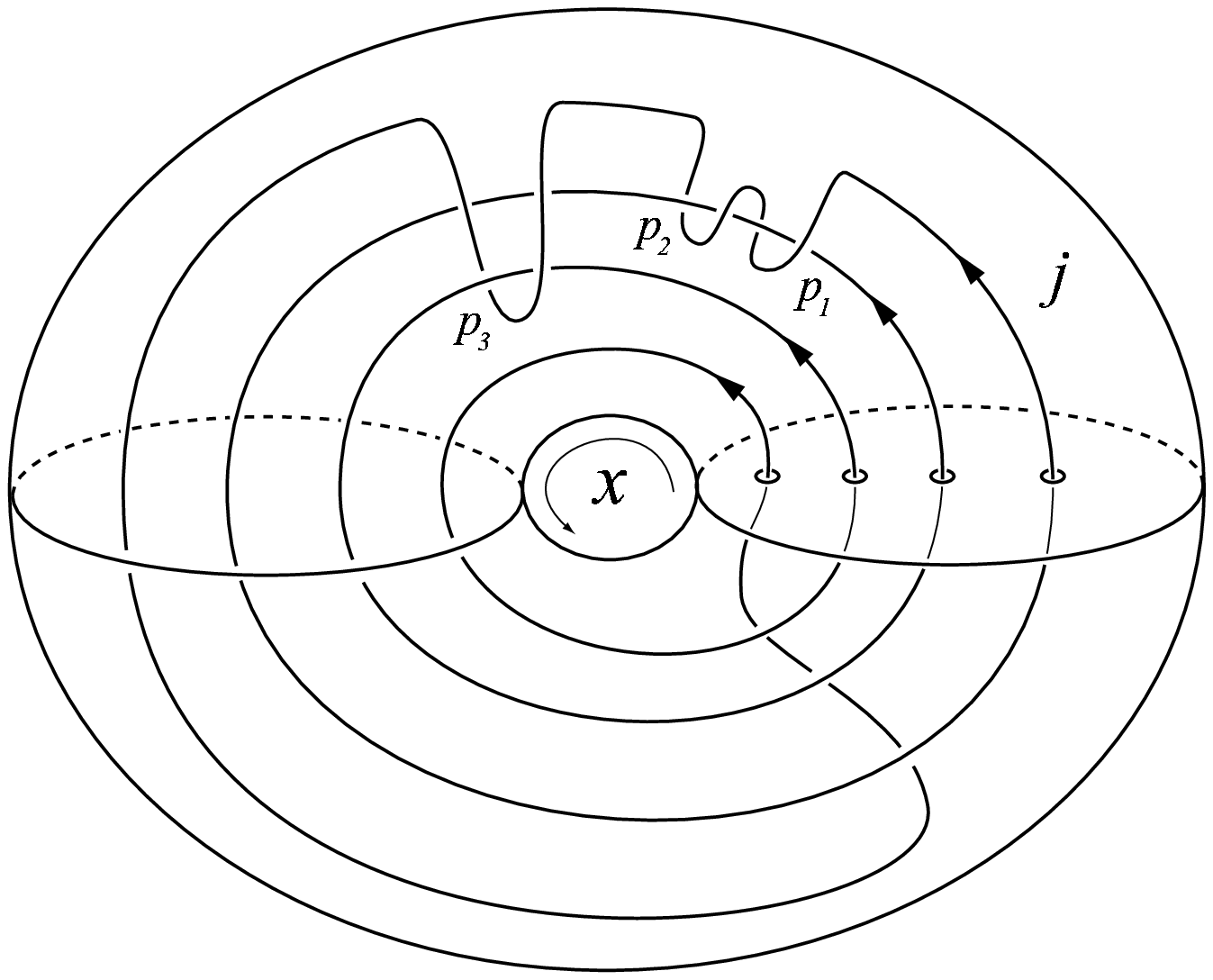}}
         \nocolon \caption{}
         \label{S1crossS2-B-fig}

\end{figure}

For $n=4$, the relation $x^2=-2x$ reduces the target space to $\Z$ (generated by the coefficient of
$x$). For example, the knot $j$ pictured in Figure~\ref{S1crossS2-B-fig} is gotten from the knot
$k$ of Figure~\ref{S1crossS2-fig} by three crossing changes (labelled $p_1$, $p_2$ and $p_3$) of
the same sign which describe the trace $H$ of a homotopy from $k$ to $j$ with three singularities.
Computing the corresponding group elements yields $\mu_k(j)=\mu(H)=x+x+x^2$ which vanishes in the
target space. We conclude that $k$ and $j$ are $W\!$--equivalent. A singular concordance admitting
Whitney disks can be constructed by tubing an $S^2$ into $H$ and performing a cusp homotopy.

\subsection{Non-spherical
example}\label{not-spherical-knot-example}

\begin{figure}[ht!]
         \centerline{\includegraphics[scale=0.55]{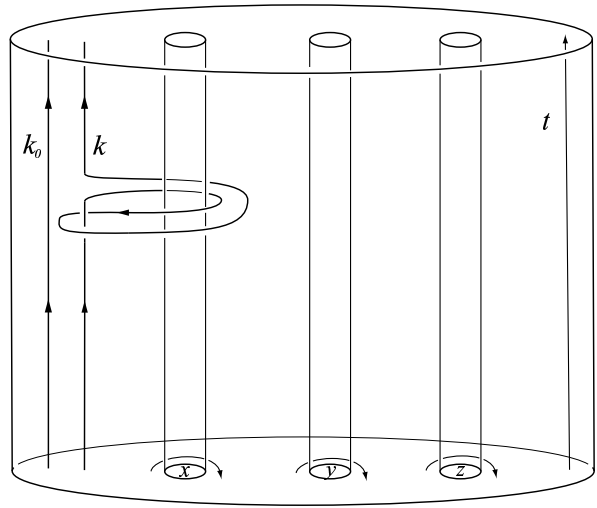}}
         \caption{A spherical knot $k_0$ and a non-spherical knot $k$}
         \label{non-spherical-knot-K-fig}

\end{figure}

\begin{figure}[ht!]
         \centerline{\includegraphics[scale=0.67]{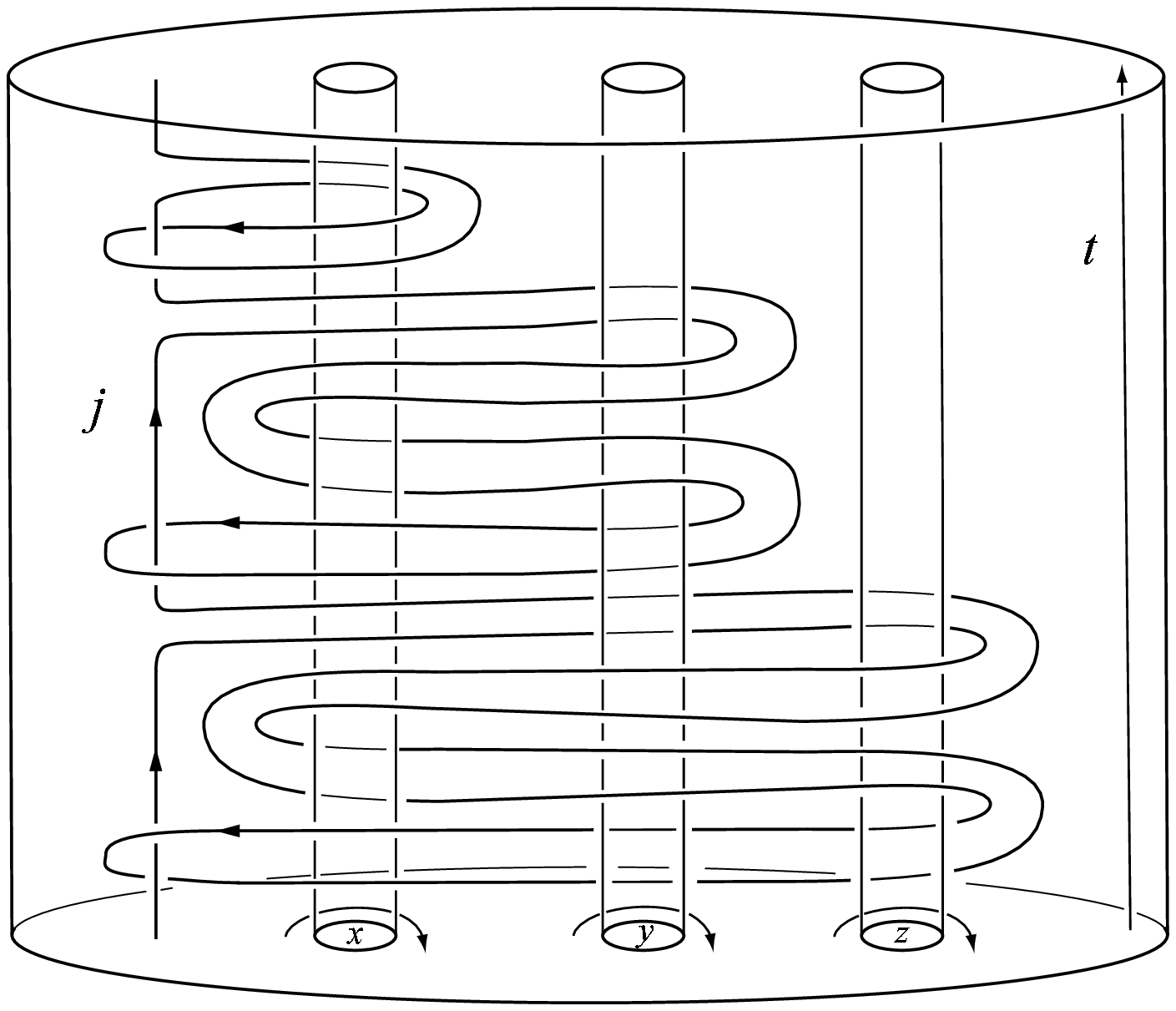}}
         \nocolon \caption{}
         \label{non-spherical-knot-J-fig}

\end{figure}

Illustrated in Figure~\ref{non-spherical-knot-K-fig} are two knots $k_0$ and $k$ in
$\mathcal{K}_t(M)$ where $M\in\mathcal{M}$ is the product $F\times S^1$ of a thrice punctured disk
$F$ with the circle and $t\in\pi_1(M)=\langle x,y,z \rangle\times\langle t\rangle$ is represented
by a circle fiber. As described in the proof of Theorem~\ref{M-thm}, $k_0$ can be self-isotoped
around each generator of $\zeta(t)=\pi_1M$ so $k_0$ is spherical. Since $M$ is irreducible, the
action of $\Phi(k_0)$ reduces to conjugation (by Proposition~\ref{Phi-prop}) and the target space
for relative self-linking numbers with respect to $k_0$ is the orbit space
$\widetilde{\Lambda}_t/({\pi_1M})$ of $\widetilde{\Lambda}_t$ under conjugation by
$\zeta(t)=\pi_1M$. Computing with the trace $H_0$ of the obvious homotopy (through one clasp
singularity) from $k_0$ to $k$ yields $\mu_{k_0}(k)=x\neq 0\in\widetilde{\Lambda}_t/({\pi_1M})$
(for some fixed orientation convention). If $K_0$ is the trace of any self-isotopy of $k_0$ with a
latitude determining $\alpha\neq 1$, then $K_{\alpha}:=-H_0+K_0+H_0$ is the trace of a
self-homotopy of $k$ having a latitude that determines $\alpha$ such that
$$
\mu(K_{\alpha})=-x+\alpha x\alpha^{-1}\neq
0\in\widetilde{\Lambda}_t. \eqno(\mbox{**})
$$
(Here we are using the same whisker for $k_0$ and $k$ as well as
assuming that latitudes of $H_0$ are trivial.) We conclude that
$k$ is {\em not} spherical (by Lemma~\ref{Phi-lem} and the
irreducibility of $M$).

Consider now the knot $j\in\mathcal{K}_t(M)$ pictured in Figure~\ref{non-spherical-knot-J-fig}.
Using the trace $J_0$ of the obvious homotopy (through 3 clasp singularities) from $k_0$ to $j$,
one computes $\mu_{k_0}(j)=x+yxy^{-1}-yzx(yz)^{-1}\in\widetilde{\Lambda}_t/({\pi_1M})$. By
Theorem~\ref{spherical-thm}, $j$ is not $W\!$--equivalent to $k$ since
$\mu_{k_0}(j)\neq\mu_{k_0}(k)\in\widetilde{\Lambda}_t/({\pi_1M})$, that is, since
$\mu(J_0)=x+yxy^{-1}-yzx(yz)^{-1}\in\widetilde{\Lambda}_t$ is not conjugate to
$\mu(H_0)=x\in\widetilde{\Lambda}_t$.

Computing with the trace $H$ of the obvious homotopy (through 2
clasp singularities) from $k$ to $j$ yields
$$
  \begin{array}{rcl}
  \mu_k(j)=\mu(H) & = & -yxy^{-1}+(yz)x(yz)^{-1} \\
   & = & y(-x+zxz^{-1})y^{-1}\\
   & = & y((\mu(K_z),z)\cdot 0)y^{-1}\\
   & = & 0\in\widetilde{\Lambda}_t/({\pi_1M}).\\
\end{array}
$$
where $K_z$ is the trace of a self-homotopy of $k$ as in $(^{**})$ above. Thus, $\mu_k(j)$ vanishes
even though $k$ and $j$ are not $W\!$--equivalent, illustrating the necessity of spherical knots in
the hypotheses of Theorem~\ref{spherical-thm}.

\subsection{A family of 3-component links}\label{3-comp-links}
Here we consider knots in the complement $M$ of a 2--component unlink $k_1\cup k_2$ in $S^3$, that
is, $M=S^1\times D^2\sharp S^1\times D^2\in\mathcal{M}$ with $\pi_1M=\langle x,y\rangle$ free on
two generators represented by meridians to $k_1\cup k_2$. Any splitting 2--sphere for the unlink
represents a generator $\sigma$ of $\pi_2M$ as a module over $\langle x,y\rangle$.

If $\gamma=\prod_{i=1}^{r}x^{m_i}y^{n_i}\in\pi_1M$, then, for any $k\in\mathcal{K}_{\gamma}(M)$,
$\widetilde{\lambda}(\sigma,k)$ is equal to
$$
\begin{array}{c}
  1-(x^{m_1})^{-1}+(x^{m_1}y^{n_1})^{-1}-(x^{m_1}y^{n_1}x^{m_2})^{-1}
+\ldots \\
  \ldots+(x^{m_1}y^{n_1}\cdots
x^{m_{r-1}}y^{n_{r-1}})^{-1}-(x^{m_1}y^{n_1}\cdots x^{m_r})^{-1}\in\widetilde{\Lambda}_{\gamma}.
\end{array}
$$

Since $\pi_1M$ is free, there exist spherical knots in every homotopy class since $\zeta(\gamma)$
is cyclic for all $\gamma\in\pi_1M$ (Remark~\ref{longitudinal-rem}). In the case that $\gamma$ is
primitive, then every knot $k$ in $\mathcal{K}_{\gamma}(M)$ is spherical and the conjugation action
of $\zeta(\gamma)=\langle\gamma\rangle$ on $\widetilde{\Lambda}_{\gamma}$ is trivial so relative
self-linking with respect to $k$ maps $\mathcal{C}_{\gamma}(M)$ (and $W\!$--equivalence classes)
onto the group generated by $\widetilde{\Lambda}_{\gamma}$ modulo the infinite family of relations
$\widetilde{\lambda}(g\cdot\sigma,k)=g\widetilde{\lambda}(\sigma,k)=0$, where $g$ ranges over
$\pi_1M$.

For example, if $\gamma=xy$, then we have the relations
$g\widetilde{\lambda}(\sigma,k)=g-gx^{-1}=g-gy=0\in\widetilde{\Lambda}_{xy}$, and one can check
(for instance, by induction on the length of reduced double coset representatives) that in this
case $\widetilde{\Lambda}_{xy}/\Phi(k)=0$ so that {\em any} two knots in $\mathcal{K}_{xy}(M)$ are
$W\!$--equivalent.

\begin{figure}[ht!]
         \centerline{\includegraphics[scale=0.50]{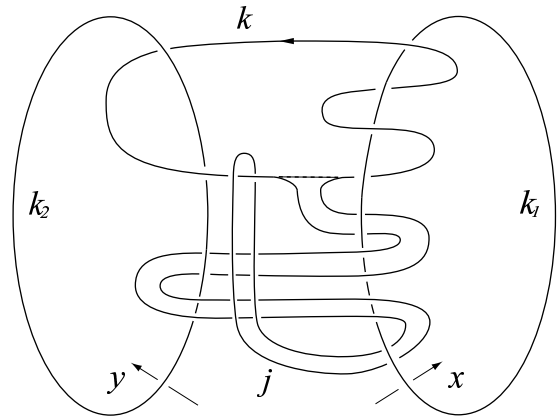}}
         \caption{The dotted arc corresponds to the knot $k$ and the ``finger'' corresponds to the knot $j$.}
         \label{xy-fig}

\end{figure}
The case $\gamma=x^2y$ shows that relative self-linking numbers are not in general determined by
projecting to the irreducible factors of $M$: Here the relations are
$g-gx^{-2}=g-gy=0\in\widetilde{\Lambda}_{x^2y}/\Phi(k)$. One can check that $xyx^{-1}$ is non-zero
in $\widetilde{\Lambda}_{x^2y}/\Phi(k)$ (which is infinitely generated) and pictured in
Figure~\ref{xy-fig} are knots $k$ and $j$ in $\mathcal{K}_{x^2y}(M)$ with $\mu_k(j)=xyx^{-1}\neq
0$. We conclude that $k$ and $j$ are not $W\!$--equivalent, a fact that is not detected by the
quotient relative self-linking numbers that ignore either one of the components of $k_1\cup k_2$
and which take values in the trivial group and $\Z_2$ respectively. (In fact, the links $k\cup k_i$
and $j\cup k_i$ are isotopic in $S^3$ for $i=1,2$.) This illustrates the ``non-abelian'' nature of
relative self-linking numbers.

A more complex example is given by $\gamma=[x,y]:=xyx^{-1}y^{-1}$ (which includes the case where
$k$ is one component of the Borromean rings). Here the target group for relative self-linking
numbers with respect to any $k\in\mathcal{K}_{[x,y]}(M)$ is described by the relations
$$
g-gx^{-1}+gx^{-1}y^{-1}-gy^{-1}=0\in\widetilde{\Lambda}_{[x,y]}/\Phi(k)
$$
where $g$ ranges over $\langle x,y\rangle$.

One can check that
$$
xy=y^{-1}x^{-1}=yx=x^{-1}y^{-1}=xy^{-1}=yx^{-1}=x+y
$$
in $\widetilde{\Lambda}_{[x,y]}/\Phi(k)$ which appears to ``almost'' split apart all cross terms;
however, I conjecture that
$$
x^{-1}y=y^{-1}x\neq x\pm y
$$
in $\widetilde{\Lambda}_{[x,y]}/\Phi(k)$ and that some cross term generators cannot be eliminated.
(A ``reason'' for this might have to do with the fact that $x^{-1}y$ and $y^{-1}x$ do not occur as
sub-words of $\gamma$ or $\gamma^{-1}$.) It would be nice to have algebraic techniques to handle
what appears to be a subtly non-abelian case.

It is worth noting that projecting to the {\em free Milnor group} (which adds the relations
$[x,x^g]=1=[y,y^g]$, $g\in\langle x,y\rangle$) centralizes $[x,y]$, and the resulting quotient of
relative self-linking numbers is indeed determined by its (abelian) projections to the irreducible
sub-manifolds and its vanishing implies that the 3-component links $k\cup k_1\cup k_2$ and $j\cup
k_1\cup k_2$ are in fact {\em disjointly $W\!$--equivalent}, meaning that their components co-bound
in $S^3 \times I$ disjoint properly immersed annuli admitting Whitney disks. It is clear that
relative self-linking numbers give isotopy invariants of classical links; a study of the more
subtle question of connections with classical link concordance invariants would be interesting.

\begin{rem}\label{computing-rem} The examples of~\ref{3-comp-links} illustrate that computing in the
infinitely presented groups $\widetilde{\Lambda}_{\gamma}/\Phi(k)$ is an interesting algebraic problem from a
topological point of view: On the one hand, the data is given purely in algebraic terms and on the other hand,
concordance classes of knots map onto the groups which are generated by double coset spaces which themselves
correspond to intersection invariants of non-simply-connected manifolds.
\end{rem}

\section{2--component
links}\label{link-sec} This section describes absolute and relative linking numbers for
2--component links which are straightforward analogues of the self-linking and relative
self-linking numbers previously defined for knots. A generalized intersection invariant $\lambda$
for immersed annuli plays the role of the generalized Wall self-intersection invariant $\mu$ and
indeterminacies are factored out in exactly the same way as for knots with analogous results in
many cases. Since the techniques and arguments involved are so similar, most proofs and details
will be omitted.

For links of null-homotopic knots, the (absolute) algebraic linking number $\lambda$ will play a
role in defining a connected sum operation for null-homotopic knots  that yields a group structure
on certain equivalence classes of null-homotopic knots (\ref{group-cor}).

Two other points worth noting are the geometric characterization of relative linking with respect to a spherical
link in terms of {\em disjointness} (\ref{spherical-link-prop}) and the modifications of Theorem~\ref{M-thm}
needed to formulate an analogous theorem for links (Theorem~\ref{link-M-thm}) as illustrated in
Example~\ref{non-spherical-link-example}.

\begin{defi}
Let $\mathcal{L}_{\gamma,\delta}(M)$ denote the set of two component links $k_1\cup k_2$ with $k_1\in
\mathcal{K}_{\gamma}(M)$ and $k_2\in \mathcal{K}_{\delta}(M)$. Let $\mathcal{C}_{\gamma,\delta}(M)$ denote
$\mathcal{L}_{\gamma,\delta}(M)$ modulo link concordance.
\end{defi}

\subsection{Links of null-homotopic knots}
We first consider 2--component links of null-homotopic knots.

\begin{defi}\label{absolute-link-defi}
For any $M$ and any link $k_1\cup k_2\in \mathcal{L}_{1,1}(M)$, define the {\em (algebraic) linking
number} $\lambda(k_1\cup k_2)$ of $k_1\cup k_2$ by
$$
\lambda(k_1\cup k_2):=\lambda(D_1,D_2)\in \Lambda=\Z[\pi_1M]
$$
where $\lambda(D_1,D_2)$ is the intersection number (\ref{wall-form}) of any properly immersed
2--disks $D_1$ and $D_2$ in $M\times I$ bounded by $k_1\cup k_2\subset M\times \{0\}$.

\end{defi}

The proof of the following theorem follows directly from the arguments in Section~\ref{null-homotopic-sec}.
\begin{thm}\label{absolute-link-thm}

For any oriented 3--manifold $M$, algebraic linking numbers induce a well-defined map from
$\mathcal{C}_{1,1}(M)$ onto $\Lambda$.
\endproof
\end{thm}

Since absolute linking invariants contain no new indeterminacies, their geometric properties are
inherited directly from Wall's intersection invariant (\ref{wall-form}); as a corollary we get that
absolute linking and self-linking numbers characterize {\em disjoint} $W\!$--equivalence to the
unlink:

\begin{cor}\label{null-geo-cor}
For any $M$ and any link $k_1\cup k_2\in \mathcal{L}_{1,1}(M)$, the following are equivalent:
\begin{enumerate}
\item The self-linking numbers $\mu(k_1)$, $\mu(k_2)$ and the linking number
$\lambda(k_1\cup k_2)$ all vanish.
\item There exist disjoint properly immersed disks $D_1$ and $D_2$
in $M\times I$ bounded by $k_1\cup k_2$ in $M\times \{0\}$ with Whitney disks pairing all self-intersections of
$D_1$ and $D_2$.
\end{enumerate}
\end{cor}
\begin{proof} Corollary~\ref{null-geo-cor} follows directly from
Proposition~\ref{wall-int-prop}. The comment before the statement of Corollary~\ref{null-geo-cor}
regarding disjoint $W\!$--equivalence to the unlink is confirmed by observing that given properly
immersed disks in $M\times I$ bounded by a link in $M\times \{0\}$, one can construct a singular
concordance to the unlink by removing a small open disk from each immersed disk followed by an
isotopy (rel $\partial$).
\end{proof}

\subsubsection{An equivalence relation on $\mathbf{\mathcal{K}_1(M)}$
admitting a group structure}

The notion of $W\!$--equivalence corresponds to the singularities of a singular concordance coming
in cancelling pairs with null-homotopic Whitney circles (\cite{FQ}). If we weaken this notion by
only requiring that the Whitney circles bound tori (rather than disks) then the resulting
equivalence relation (for knots) is characterized by a slightly weaker absolute self-linking number
$\mu^{\pi}(k)$ which takes values in the free abelian group $\widetilde{\Lambda}^{\pi}$ generated
by the non-trivial conjugacy classes of $\pi_1M$ modulo inversion (a quotient of
$\widetilde{\Lambda}$). The resulting equivalence classes $\mathcal{K}^{\pi}_1(M)$ in
$\mathcal{K}_1(M)$ can be made into a group:

\begin{cor}\label{group-cor}
There exists a connected sum operation $\sharp$ that makes $\mathcal{K}_1^{\pi}(M)$ into a group isomorphic to
$\widetilde{\Lambda}^{\pi}$, the isomorphism being given by $\mu^{\pi}$.
\end{cor}

\begin{proof} The connected sum operation $\sharp$ of
Corollary~\ref{group-cor} will be defined after giving the precise definition of $\mu^{\pi}(k)$.
The usual definition of connected sum for knots in $S^3$ (which makes $\mathcal{C}(S^3)$ into an
abelian group) involves two steps; first the knots are separated by any embedded 2--sphere, then
they are joined by the (orientation-preserving) band sum (0-surgery) along any choice of band that
meets the 2--sphere in a standard way. Here we first achieve an ``algebraic separation'' by
requiring the vanishing of the absolute algebraic linking number, then join the knots along any
generic band.

Define $\widetilde{\Lambda}^\pi$ by
$$
\widetilde{\Lambda}^\pi:=\frac{\Z[\pi_1M]}{\Z[1]\oplus\langle g-xg^{\pm 1}x^{-1}\rangle};\quad x,g\in\pi_1M
$$
that is, $\widetilde{\Lambda}^{\pi}$ is the free abelian group generated by the non-trivial conjugacy classes of
$\pi_1M$ modulo inversion. ($\widetilde{\Lambda}^{\pi}$ is the quotient of $\widetilde{\Lambda}$ by the relations
$g-xgx^{-1}=0$ for $x,g\in\pi_1M$.)

\begin{defi}
For $k\in\mathcal{K}_1(M)$, define
$$
\mu^{\pi}(k):=\mu(D)\in\widetilde{\Lambda}^{\pi}
$$
where $D$ is any singular null-concordance for $k$.

If $\mu^{\pi}(k)=\mu^{\pi}(j)$ we say that $k$ and $j$ are {\em $\mu^{\pi}\!$--equivalent}. Denote
by $\mathcal{K}_1^{\pi}(M)$ the set $\mathcal{K}_1(M)$ modulo $\mu^{\pi}\!$--equivalence.
\end{defi}

It follows from Theorem~\ref{spherical-thm} (and the fact that the unknot is spherical) that two
knots are $\mu^{\pi}\!$--equivalent if and only if they co-bound a singular concordance with all
singularities paired by ``Whitney tori'' (rather than disks) since the singularities will occur in
cancelling pairs with null-homologous Whitney-circles (rather than null-homotopic) due to the
relations $g-xgx^{-1}=0$ added in the definition of $\mu^{\pi}(k)$.

For knots $k_1$ and $k_2$ in $\mathcal{K}_1^{\pi}(M)$ define $k_1\sharp k_2$ as follows:

First perturb $k_1$ to be disjoint from $k_2$ in $M$, then isotope $k_1$ so that
$\lambda(k_1,k_2)=0\in{\Lambda}$. This isotopy can be assumed to be supported in the neighborhood
of a finite collection of arcs which guide crossing changes between $k_1$ and $k_2$. Of course this
isotopy does not change the class of $k_1$ or $k_2$ in $\mathcal{K}_1^{\pi}(M)$ (or in
$\mathcal{K}_1(M)$).

Now define $k_1 \sharp k_2$ to be the usual orientation preserving band sum (0-surgery) of $k_1$ with $k_2$ along
any band $b$ (whose interior is disjoint from $k_1\cup k_2$).

If $D_1$ and $D_2$ are any singular null-concordances of $k_1$ and $k_2$, then the band sum
$D_1\sharp_bD_2$ is a singular null-concordance for $k_1\sharp k_2$. Taking the whisker on $D_1$ as
a whisker for $D_1\sharp_bD_2$ we have
$$
\mu(D_1\sharp_bD_2)=\mu(D_1)+\beta\mu(D_2)\beta^{-1}+\lambda(k_1,k_2)\beta^{-1}\in\widetilde{\Lambda}
$$
where $\beta\in\pi_1M$ is determined by the band $b$ and the whiskers on $k_1$ and $k_2$. Since $\lambda(k_1,k_2)$
vanishes, taking this equation in $\widetilde{\Lambda}^{\pi}$ yields
$$
\mu^{\pi}(k_1\sharp k_2)=\mu^{\pi}(k_1)+\mu^{\pi}(k_2)
$$
which shows that $\sharp$ is well-defined on $\mathcal{K}_1^{\pi}(M)$. Inverses are gotten by
reversing knot orientations; $\mu^{\pi}$ is injective by definition and surjective by
Theorem~\ref{absolute-link-thm}.
\end{proof}

\subsection{Links of essential knots}\label{essential-link-subsec}
The definition of relative linking numbers for links of essential knots follows directly along the lines of the
case of relative self-intersection numbers for knots. The rest of this section contains statements of the
analogous theorems without proofs.

\subsection{Intersection numbers for annuli}\label{annuli-ints}
Let $A$ and $B$ be a pair of properly immersed annuli in a 4--manifold $X$. Chose whiskers for each
of $A$ and $B$ and let $\langle \gamma \rangle$ and $\langle \delta \rangle$ be the respective
cyclic images in $\pi_1X$ of the induced maps on fundamental groups. The group element associated
to an intersection point from a loop through $A$ and back through $B$ in the usual way is only
well-defined up to left multiplication by powers of $\gamma$ and right multiplication by powers of
$\delta$. Thus, we work with the double cosets of $\pi_1X$ by $\langle \gamma \rangle$ on the left
and $\langle \delta \rangle$ on the right. Adjusting the notation of ~\ref{annuli-self-ints}, we
define

$$
\Lambda_{\gamma,\delta}:=\Z[{\langle \gamma \rangle}\backslash
\pi_1X/{\langle \delta \rangle}] =\frac{\Z[\pi_1X]}{\{g-{\gamma}^n
g {\delta}^m\}}
$$
where $n$ and $m$ range over the integers.

\begin{defi}
Let $A$ and $B$ be properly immersed annuli (equipped with whiskers) in a 4--manifold $X$ whose
fundamental groups map to $\langle \gamma \rangle$ and $\langle \delta \rangle$ respectively. The
{\em intersection number} $\lambda(A,B)$ of $A$ and $B$ is defined by
$$
\lambda(A,B):=\sum (\sign  p)\cdot g_p \in \Lambda_{\gamma,\delta}
$$
where the sum is over all intersection points $p\in A \cap B$ and
$\sign  p$ comes from the orientations of $X$ and the sheets of
$A$ and $B$ at $p$ as usual.

\end{defi}
Note that $\mu(A,B)$ is well-defined up to left (resp.\ right) multiplication by elements of
$\zeta(\gamma)$ (resp.\ $\zeta(\delta)$) corresponding to the effect of changing the whisker of $A$
(resp.\ $B$) while preserving the image of the induced map on fundamental groups.

A whisker change that does not preserve the cyclic images also
changes the target space by an isomorphism.

\begin{prop}\label{annuli-int=0}

The intersection number $\lambda(A,B)\in \Lambda_{\gamma,\delta}$
is invariant under homotopy (rel $\partial$) and vanishes if and
only if there exists Whitney disks pairing all intersections
between $A$ and $B$ if and only if $A$ and $B$ can be homotoped
(rel $\partial$) to be disjoint. \endproof
\end{prop}

\subsection{The indeterminacy subgroup $\Phi(k_1\cup k_2)$}

\begin{defi}\label{Phi-link-defi}

For any link $k_1\cup k_2\in\mathcal{L}_{\gamma,\delta}(M)$, define
$$
\Phi(k_1\cup k_2)\leq {\Lambda}_{\gamma,\delta}\rtimes
(\zeta(\gamma)\times\zeta(\delta))
$$
to be the subgroup of the semi-direct product of
${\Lambda}_{\gamma,\delta}$ and $\zeta(\gamma)\times\zeta(\delta)$
with respect to the action
$$
(\phi,\psi):y\mapsto\phi y \psi^{-1}
$$
of $\zeta(\gamma)\times\zeta(\delta)$ on
${\Lambda}_{\gamma,\delta}$ generated by the elements
$$
(\lambda(K_1,K_2)(w_1,w_2),(\lat[K_1](w_1),\lat[K_2](w_2))
$$
as $K_1\cup K_2$ ranges over all singular self-concordances of
$k_1\cup k_2$ and $w_1$ (resp.\ $w_2$) ranges over all whiskers
identifying $[k_1]=\gamma\in\pi_1M$ (resp.\
$[k_2]=\delta\in\pi_1M$).

\end{defi}

\subsubsection{The action of $\Phi(k_1\cup k_2)$ on
${\Lambda}_{\gamma,\delta}$}

The group ${\Lambda}_{\gamma,\delta}\rtimes
(\zeta(\gamma)\times\zeta(\delta))$ acts on
${\Lambda}_{\gamma,\delta}$ by
$$
(z,(\phi,\psi)): y \mapsto z+\phi y \psi^{-1}
$$
for $(z,(\phi,\psi))\in{\Lambda}_{\gamma,\delta}\rtimes
(\zeta(\gamma)\times\zeta(\delta))$ and
$y\in{\Lambda}_{\gamma,\delta}$. We denote the equivalence classes
under the restriction of this action to the subgroup $\Phi(k_1\cup
k_2)$ by
$$
{\Lambda}_{\gamma,\delta}/\Phi(k_1\cup k_2).
$$

\subsubsection{The action of $\zeta(\gamma)\times\zeta(\delta)$ on
${\Lambda}_{\gamma,\delta}/\Phi(k_1\cup k_2)$}

The effect that changing the whiskers for $k_1\cup
k_2\in\mathcal{L}_{\gamma,\delta}(M)$ has on a generator of
$\Phi(k_1\cup k_2)$ is described by an action of
$\zeta(\gamma)\times\zeta(\delta)$ on
${\Lambda}_{\gamma,\delta}/\Phi(k_1\cup k_2)$ defined by:
$$
(\alpha,\beta):(z,(\phi,\psi))\mapsto (\alpha z \beta^{-1},
(\alpha \phi\alpha^{-1},\beta\psi\beta^{-1}))
$$
for all $(\alpha,\beta)\in\zeta(\gamma)\times\zeta(\delta)$.

We have a well-defined conjugation action of
$\zeta(\gamma)\times\zeta(\delta)$ on
${\Lambda}_{\gamma,\delta}/\Phi(k_1\cup k_2)$ since if
$$
x= z+\phi y \psi^{-1} \in{\Lambda}_{\gamma,\delta}
$$
we have
$$
\alpha x\beta^{-1}=\alpha
z\beta^{-1}+(\alpha\phi\alpha^{-1})(\alpha
y\beta^{-1})(\beta\psi^{-1}\beta^{-1})\in
{\Lambda}_{\gamma,\delta}.
$$

Denote by
$$
{\Lambda}_{\gamma,\delta}/(\Phi(k_1\cup
k_2),\zeta(\gamma)\times\zeta(\delta))
$$
the orbit space of ${\Lambda}_{\gamma,\delta}/\Phi(k_1\cup k_2)$
under the action of $\zeta(\gamma)\times\zeta(\delta)$.

\subsection{Definition of relative linking numbers
$\lambda_{k_1\cup k_2}(j_1\cup j_2)$}

\begin{defi}\label{link-defi}
Fix a link $k_1\cup k_2$ in $\mathcal{L}_{\gamma,\delta}(M)$. For all $j_1\cup
j_2\in\mathcal{L}_{\gamma,\delta}(M)$, define $\lambda_{k_1\cup k_2}(j_1\cup j_2)$, the {\em
relative linking number of $j_1\cup j_2$ with respect to $k_1\cup k_2$}, by
$$
\lambda_{k_1\cup k_2}(j_1\cup j_2):=\lambda(H_1,H_2)
\in{\Lambda}_{\gamma,\delta}/(\Phi(k_1\cup
k_2),\zeta(\gamma)\times\zeta(\delta))
$$
where $H_1\cup H_2$ is any singular concordance from $k_1\cup k_2$
to $j_1\cup j_2$.
\end{defi}

\begin{thm}\label{rel-link-thm}
For each $k_1\cup k_2\in\mathcal{L}_{\gamma,\delta}(M)$, the map $j_1\cup j_2\mapsto
\lambda_{k_1\cup k_2}(j_1\cup j_2)$ induces a well-defined map
$$
\mathcal{C}_{\gamma,\delta}(M)\twoheadrightarrow{\Lambda}_{\gamma,\delta}/(\Phi(k_1\cup
k_2),\zeta(\gamma)\times\zeta(\delta))
$$
onto the target space.\endproof
\end{thm}

\subsection{A generating set for $\Phi(k_1\cup k_2)$}
\begin{prop}\label{Phi-link-prop}

Given $k_1\cup k_2\in\mathcal{L}_{\gamma,\delta}(M)$, let
$\{\phi_q\}$ (resp.\ $\{\psi_r\}$) be any generating set for the
centralizer $\zeta(\gamma)$ (resp.\ $\zeta(\delta)$) and let
$\{K^1_q\}$ (resp.\ $\{K^2_r\}$) be a corresponding set of traces
of self-homotopies of $k_1$ (resp.\ $k_2$) such that, for a fixed
whisker $w_1$ on $k_1$ (resp.\ $w_2$ on $k_2$),
$\lat[K^1_q](w_1)=\phi_q$ for all $q$ (resp.\
$\lat[K^2_r](w_2)=\psi_r$ for all $r$). Denote the traces
$k_1\times I$ and $k_2\times I$ of constant self-homotopies of
$k_1$ and $k_2$ by $K^1$ and $K^2$. Let $\{\sigma_s\}$ be any
generating set for $\pi_2M$ as a module over $\pi_1M$.

Let $\mathcal{S}_0\subset{\Lambda}_{\gamma,\delta}\rtimes
(\zeta(\gamma)\times\zeta(\delta))$ denote the union of the sets
$$
\{(\lambda(K^1_q,K^2)(w_1,w_2),(\phi_q,1)),
(\lambda(K^1,K^2_r)(w_1,w_2),(1,\psi_r)),\}
$$
and
$$
\{(g\lambda(\sigma_s,k_2)(w_1,w_2),(1,1)),
(\lambda(k_1,\sigma_s)(w_1,w_2)g,(1,1)) \}
$$
where $g$ ranges over $\pi_1M$.

Then $\Phi(k_1\cup k_2)$ is generated by
$$
\bigcup_{(\alpha,\beta)\in\zeta(\gamma)\times\zeta(\delta) }\alpha\mathcal{S}_0\beta^{-1}
$$
that is, by the union of all the sets in the orbit of
$\mathcal{S}_0$ under the action of
$\zeta(\gamma)\times\zeta(\delta)$.\endproof

\end{prop}

\begin{figure}[ht!]
         \centerline{\includegraphics[scale=0.60]{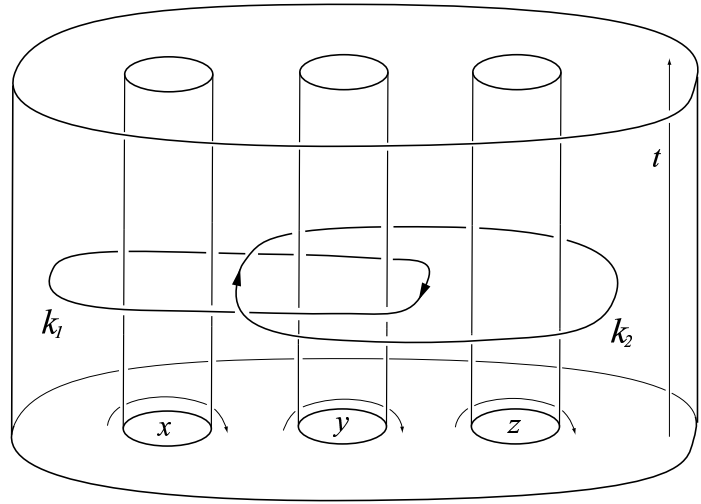}}
         \caption{A non-spherical link in the product of a thrice punctured 2--disk with the circle}
         \label{nonspherical-link-example-fig}

\end{figure}

\subsection{Spherical links}
\begin{defi}\label{spherical-link-defi}
A link $k_1\cup k_2\in\mathcal{L}_{\gamma,\delta}(M)$ is {\em
spherical} if there exist whiskers $w_1$ and $w_2$ on $k_1$ and
$k_2$ such that for any
$(\phi,\psi)\in\zeta(\gamma)\times\zeta(\delta)$ there exist
singular self-concordances $K^1_{\phi}$ of $k_1$ and $K^2_{\psi}$
of $k_2$ with
$\lambda(K^1_{\phi},K^2_{\psi})=0\in\Lambda_{\gamma,\delta}$ and
$\lat[K^1_{\phi}](w_1)=\phi$ and $\lat[K^2_{\psi}](w_2)=\psi$.
\end{defi}

\begin{prop}\label{spherical-link-prop}
If $k_1\cup k_2\in\mathcal{L}_{\gamma,\delta}(M)$ is spherical, then for any links $j_1\cup j_2$
and $j'_1\cup j'_2$ in $\mathcal{L}_{\gamma,\delta}(M)$ the following are equivalent:
\begin{enumerate}
  \item $\lambda_{k_1\cup k_2}(j_1\cup j_2)=\lambda_{k_1\cup k_2}(j'_1\cup j'_2)$.
  \item  There
exists a singular concordance $H^1\cup H^2$ between $j_1\cup j_2$
and $j'_1\cup j'_2$ such that $H^1\cap H^2=\emptyset$.\endproof
  \end{enumerate}
  \end{prop}

Spherical links are ``less plentiful'' than spherical knots: Even when restricting to manifolds in
$\mathcal{M}$ (Section~\ref{M-sec}), there may exist pairs $\gamma,\delta$ of elements in $\pi_1M$
with $\mathcal{L}_{\gamma,\delta}(M)$ containing no spherical links as described in
Example~\ref{non-spherical-link-example} below.

We do, however, have the following result which says roughly that if the centralizers of $\gamma$
and $\delta$ are either ``small enough'' or ``large enough'' then $\mathcal{L}_{\gamma,\delta}(M)$
does contain spherical links.\eject

\begin{thm}\label{link-M-thm}
For any $M$, there exists a spherical link in
$\mathcal{L}_{\gamma,\delta}(M)$ if any one of the following is
satisfied:
\begin{enumerate}
  \item Both of the centralizers $\zeta(\gamma)$ and
  $\zeta(\delta)$ are cyclic.
  \item Both of $\gamma$ and $\delta$ are carried by regular fibers
  of the characteristic Seifert fibered sub-manifolds which carry $\zeta(\gamma)$ and
  $\zeta(\delta)$.
  \item At least one of $\gamma$ and $\delta$ are equal to
  $1\in\pi_1M$.\endproof
\end{enumerate}
\end{thm}
Theorem~\ref{link-M-thm} can be proved along the lines of the
proof of Theorem~\ref{M-thm} by finding disjoint self-isotopies of
the spherical knots described there whose latitudes represent
generators of $\zeta(\gamma)\times\zeta(\delta)$.

\begin{example}\label{non-spherical-link-example}
For $M=F\times S^1$ as in Example~\ref{not-spherical-knot-example}, if $\gamma=xy$ and $\delta=yz$,
then $\mathcal{L}_{\gamma,\delta}(M)$ contains no spherical links. As illustrated in
Figure~\ref{nonspherical-link-example-fig}, one can see an obstruction measured by an intersection
number between curves in $F$ representing $\gamma$ and $\delta$ which can be used to compute the
indeterminacy subgroup (\cite{S2}).
\end{example}

\Addresses\recd

\end{document}